\documentclass[12pt]{svjour3}     

\makeatletter
\def\ps@pprintTitle{%
 \let\@oddhead\@empty
 \let\@evenhead\@empty
 \def\@oddfoot{}%
 \let\@evenfoot\@oddfoot}
\makeatother

\usepackage{tabularx} 
\usepackage{caption}  
\usepackage[export]{adjustbox}

\usepackage{placeins,subcaption,tabularx}
\usepackage{amsmath,amscd,wrapfig,amsfonts,amssymb,mathtools,mathrsfs}
\usepackage{booktabs,multirow,lscape,cool,bm}
\usepackage{url,parskip,caption}

\providecommand{\e}[1]{\ensuremath{\times 10^{#1}}}
\usepackage[top=1.10in, bottom=1.10in, left=1.10in, right=1.10in]{geometry}

\makeatletter
\g@addto@macro\normalsize{%
  \setlength\abovedisplayskip{.6em}
  \setlength\belowdisplayskip{.6em}
  \setlength\abovedisplayshortskip{.6em}
  \setlength\belowdisplayshortskip{.6em}
}

\begin{document}









\title{A solver for linear scalar ordinary differential equations whose running time is bounded independent of frequency}
\titlerunning{Frequency-independent ODE solver}
\journalname{Journal of Scientific Computing}

\author{Murdock Aubry         \and
        James Bremer
}

\institute{Murdock Aubry \at
              Department of Mathematics\\
              University of Toronto \\
              \email{murdock.aubry@mail.utoronto.ca}           
           \and
              James Bremer \at
              Department of Mathematics \\
              University of Toronto \\
              \email{bremer@math.toronto.edu}           
}

\date{\vskip -7em}

\maketitle

\begin{abstract}
When the eigenvalues of the coefficient matrix for a linear scalar ordinary differential equation
are of large magnitude, its solutions exhibit complicated behaviour, such as high-frequency oscillations,
rapid growth or rapid decay.  The cost of representing such solutions using standard techniques
grows with the magnitudes of the eigenvalues.  As a consequence, the  running times of most solvers for 
ordinary differential equations also grow with these eigenvalues.
However, a large class of  scalar ordinary differential equations with slowly-varying coefficients
admit slowly-varying phase functions that can be  represented at a cost which is bounded
independent of the magnitudes of the eigenvalues of the corresponding coefficient matrix.
Here, we introduce a numerical algorithm for constructing slowly-varying phase functions which represent the solutions
of  a linear scalar ordinary differential equation.  Our method's running time depends on the complexity of the equation's coefficients,
but is bounded independent of the magnitudes of the equation's eigenvalues.
Once the phase functions have been constructed, essentially any reasonable initial or boundary value problem 
for the scalar equation can be easily solved.  
We present the results of numerical experiments showing that, despite 
its greater generality, our algorithm is competitive with state-of-the-art
methods for  solving highly-oscillatory second  order differential equations.
We also compare our method with  Magnus-type exponential integrators and find
that our approach is  orders of magnitude faster in the high-frequency regime.

\end{abstract}

\begin{section}{Introduction}

The complexity of the solutions of an  $n^{th}$ order linear homogeneous ordinary differential equation
\begin{equation}
y^{(n)}(t) + q_{n-1}(t) y^{(n-1)}(t) + \cdots + q_1(t) y'(t) + q_0(t) y(t) = 0
\label{introduction:scalarode}
\end{equation}
increases with the magnitudes of the eigenvalues $\lambda_1(t), \ldots, \lambda_n(t)$
of the coefficient matrix
\begin{equation}
\left(\begin{array}{cccccccc}
0 & 1 & 0 & \cdots & 0 & 0 \\
0 & 0 & 1 & \cdots & 0 & 0 \\
\vdots &  &  & \ddots &  & \vdots\\
0 & 0 & 0 & \cdots & 1 & 0 \\
0 & 0 & 0 & \cdots & 0 & 1 \\
-q_0(t) & -q_1(t) & -q_2(t) & \cdots&  -q_{n-2}(t) & -q_{n-1}(t)\\
\end{array}
\right)
\label{introduction:scalarcoef}
\end{equation}
obtained from  (\ref{introduction:scalarode}) in the usual way.  
Indeed, the cost to represent such solutions over an interval $[a,b]$ using standard techniques
(e.g., polynomial or trigonometric expansions) typically grows roughly linearly with the
quantity
\begin{equation}
\Omega = \max_{i=1,\ldots,n} \int_a^b \left| \lambda_i(t) \right|\, dt,
\end{equation}
which we refer to as the frequency of (\ref{introduction:scalarode}).
We use this terminology
because, in most cases of interest, it is the imaginary parts of the eigenvalues which
are of large magnitude.  Indeed, when the real part of one or more of the $\lambda_j(t)$ is large in 
size, most initial and terminal value problems for (\ref{introduction:scalarode}) are highly
ill-conditioned and solving them numerically requires specialized techniques which exploit additional
information about the desired solution.

Although the complexity of the solutions of (\ref{introduction:scalarode}) increases with frequency,
a large class of linear scalar ordinary differential equations admit phase functions whose
cost to represent via standard techniques is bounded independent of the magnitudes of the eigenvalues of (\ref{introduction:scalarcoef}).
In fact, if $q_0,\ldots,q_{n-1}$ are slowly-varying on an interval $I$ and 
the differential equation (\ref{introduction:scalarode}) is nondegenerate there --- meaning that 
the  eigenvalues  $\lambda_1(t), \ldots, \lambda_n(t)$   are distinct for each $t \in I$ ---
then it is possible to find slowly-varying phase functions  $\psi_1,\ldots,\psi_n\colon I \to \mathbb{C}$ such that 
\begin{equation}
\left\{ \exp\left(\psi_j(t)\right)  :  j=1,\ldots,n \right\}
\label{introduction:phaserep}
\end{equation}
is a basis for the space of solutions of (\ref{introduction:scalarode}) given on the interval $I$.
That slowly-varying phase functions exist under these conditions, at least in an asymptotic sense,
has long been known.  Indeed, this observation is the basis of the WKB method and other related techniques 
(see, for instance, \cite{Miller}, \cite{Wasov} and \cite{SpiglerPhase1,SpiglerPhase2,SpiglerZeros}).
A theorem which establishes  the existence of  slowly-varying phase functions for 
second order differential equations under mild conditions on their coefficients is proven in \cite{BremerRokhlin}.
Although it is not immediately obvious how to generalize the argument of \cite{BremerRokhlin} to higher
order scalar equations, known results regarding the asymptotic approximation 
of solutions of differential equations and numerical evidence (including the experiments of this paper)
strongly suggest  the situation for higher order scalar equations is much the same as it is for
second order equations.

The derivatives of the phase functions $\psi_1,\ldots,\psi_n$, which we 
denote by $r_1,\ldots,r_n$, satisfy  an $(n-1)^{st}$ order  nonlinear inhomogeneous ordinary differential equation,
the general form of which is quite complicated.   When $n=2$, it is the Riccati equation
\begin{equation}
r'(t) + (r(t))^2 + q_1(t) r(t) + q_0(t) = 0;
\label{introduction:riccati1}
\end{equation}
when $n=3$, the nonlinear equation is
\begin{equation}
r''(t) + 3 r'(t) r(t) + (r(t))^3 + q_2(t) r'(t) + q_2(t) (r(t))^2 + q_1(t) r(t) + q_0(t) = 0;
\label{riccati2}
\end{equation}
and, for $n=4$, we have
\begin{equation}
\begin{aligned}
r'''(t) &+ 4 r''(t) r(t) + 3 (r'(t))^2 + 6  r'(t) (r(t))^2 + (r(t))^4 +q_3(t) (r(t))^3 + q_3(t) r''(t) \\
&+ 3 q_3(t) r'(t) r(t) +q_2(t) (r(t))^2 +q_2(t) r'(t) + q_1(t) r(t) +q_0(t) = 0.
\end{aligned}
\label{riccati3}
\end{equation}
By a slight abuse of terminology, we will refer to the $(n-1)^{st}$ order nonlinear
equation obtained by inserting the representation 
\begin{equation}
y(t) = \exp\left(\int r(t)\, dt\right)
\end{equation}
into (\ref{introduction:scalarode}) as the $(n-1)^{st}$ order Riccati equation, or, alternatively, 
the Riccati equation for (\ref{introduction:scalarode}).

An obvious approach to initial and boundary value boundary problems for (\ref{introduction:scalarode}) calls for constructing
a suitable collection of slowly-varying phase functions by solving the corresponding Riccati equation numerically.
Doing so is not as straightforward as it sounds, however.  The principal difficulty is that 
most solutions of the Riccati equation for (\ref{introduction:scalarode}) are rapidly-varying when 
the eigenvalues $\lambda_1(t),\ldots,\lambda_n(t)$ are of large magnitude, and some mechanism is needed to 
select the slowly-varying solutions.

The article \cite{BremerPhase} introduces an algorithm for constructing two slowly-varying phase function 
$\psi_1$ and $\psi_2$ such that $\exp(\psi_1(t))$ and $\exp(\psi_2(t))$ 
constitute a basis in the space of solutions of a second order linear ordinary differential equation
of the form
\begin{equation}
y''(t) + q(t) y(t) = 0,\ \ \ \ \ a <t <b,
\label{introduction:two}
\end{equation}
where $q$ is slowly-varying and non-vanishing on $(a,b)$.  It operates  by  
constructing a smoothly deformed version of the coefficient $q$ which is equal
to an appropriately chosen constant in a neighborhood of some point $c$ in $(a,b)$ and coincides with the original 
coefficient $q$ in a neighborhood of a point $d$ in $(a,b)$.  There is a pair of slowly-varying
phase functions for the deformed equation whose derivatives at $c$ are known and whose derivatives at $d$ 
coincide with the derivatives of a pair of slowly-varying phase functions for the original equation.
Consequently, by solving the Riccati equation corresponding to the deformed equation 
with initial conditions specified at $c$, the values of the derivatives of a pair of slowly-varying
phase functions for the original equation at the point $d$ can be calculated.
Once this has been done, the Riccati equation corresponding to the original equation
is solved using the values at $d$ as initial conditions
in order to calculate the derivatives of a pair of slowly-varying phase functions
for (\ref{introduction:two})   over the whole interval.  The desired slowly-varying
phase functions $\psi_1$ and $\psi_2$ are obtained by integration.  
The cost of the entire procedure is bounded independent of the magnitude of $q$, which
is related to the eigenvalues of the coefficient matrix corresponding to (\ref{introduction:two}) 
via 
\begin{equation}
\lambda_1(t) =  \sqrt{-q(t)} \ \ \ \mbox{and} \ \  \ \lambda_2(t) =  -\sqrt{-q(t)}.
\label{introduction:eigentwo}
\end{equation}

From (\ref{introduction:eigentwo}), it follows that the assumption that $q$ is non-vanishing on $(a,b)$ is equivalent to the 
condition that (\ref{introduction:two}) is nondegenerate on $(a,b)$.  
In \cite{BremerPhase2}, the  method of \cite{BremerPhase} is extended 
to the case in which  (\ref{introduction:two}) is nondegenerate on an interval $[a,b]$ except at a finite number 
of turning points.  The equation (\ref{introduction:scalarode}) has a turning point at $t_0$ 
provided the eigenvalues  $\lambda_1(t), \ldots, \lambda_n(t)$ of (\ref{introduction:scalarcoef})
are distinct in a deleted neighborhood of $t_0$, but coalesce at $t_0$. 
The turning points of (\ref{introduction:two}), then,  are precisely the isolated zeros of $q$.
Because slowly-varying phase functions need not extend across turning points, the algorithm
of \cite{BremerPhase2}   introduces a  partition $a=\xi_1 < \xi_2 < \ldots < \xi_k=b$ of 
$[a,b]$ such that $\xi_2,\ldots,\xi_{k-1}$ are the roots of $q$ in the open
interval $(a,b)$.  
It then  applies a variant of the method of \cite{BremerPhase}
to each of the subintervals $\left[\xi_j,\xi_{j+1}\right]$, $j=1,\ldots,k-1$,
which results in a collection
of $2(k-1)$ slowly-varying phase functions that efficiently represent the solutions
of (\ref{introduction:two}).  

It is relatively straightforward to generalize the approach of \cite{BremerPhase} to 
the case of nondegenerate higher order scalar equations.   However, while the resulting 
algorithm  is highly-effective for a large class of equations of the form 
(\ref{introduction:scalarode}),  the authors have found another approach 
inspired by the classical Levin method for evaluating oscillatory integrals to be somewhat more robust.
Introduced in  \cite{Levin}, the Levin method is based on the observation that
if $p_0$ and $f$ are slowly varying, then the  inhomogeneous equation
\begin{equation}
y'(t) + p_0(t) y(t)  = f(t)
\label{introduction:levin}
\end{equation}
has a slowly-varying solution $y_0$, regardless of  the magnitude of $p_0$.
Similarly to the case of phase functions,
the proofs appearing in \cite{Levin} and subsequent works on the Levin method
do not immediately apply to the case of higher order scalar equations, but 
experimental evidence and results for special cases strongly suggest that the Levin 
principle generalizes.  That is to say, equations of the form
\begin{equation}
y^{(n)}(t) + p_{n-1}(t) y^{(n-1)}(t) + \cdots + p_1(t) y'(t) + p_0(t) y(t) = f(t).
\label{introduction:inhom}
\end{equation}
admit solutions whose complexity depends on that of the right-hand side $f$ and of
the coefficients $p_0,\ldots,p_{n-1}$, but is bounded independent of the magnitudes of $p_0,\ldots,p_{n-1}$.

The algorithm of this paper exploits the existence of slowly-varying phase functions
and the Levin principle to solve initial and boundary value problems
for nondegenerate scalar equations of the form (\ref{introduction:scalarode}) with slowly-varying
coefficients. 
It  operates by constructing  slowly-varying phase functions $\psi_1\ldots,\psi_n$ such that (\ref{introduction:phaserep}) is a basis in the space
of solutions of a nondegenerate scalar equation.  Once this has been done, any reasonable
initial or boundary value problem for (\ref{introduction:scalarode}) can be solved more-or-less
instantaneously.     As with \cite{BremerPhase}, the method of this paper
can be extended to the case of a scalar equation which is nondegenerate on an interval
$[a,b]$ except at a finite number of turning points by applying it on a collection
of subintervals of $[a,b]$; however, for the sake of simplicity, we consider only
nondegenerate equations here.

The algorithms of \cite{BremerPhase}, \cite{BremerPhase2} and this article
bear some superficial similarities to Magnus expansion methods.
Introduced in \cite{Magnus}, Magnus expansions are certain series of the form
\begin{equation}
\sum_{k=1}^\infty \Omega_k(t)
\label{introduction:Magnus}
\end{equation}
such that  $\exp\left(\sum_{k=1}^\infty \Omega_k(t)\right)$
locally represents a fundamental matrix for a system of  differential equations
\begin{equation}
\mathbf{y}'(t) = A(t) \mathbf{y}(t).
\label{introduction:system}
\end{equation}
The first few terms for the series around $t=0$ are given by 
\begin{equation}
\begin{aligned}
\Omega_1(t) &=  \int_{0}^t A(s)\ ds,\\
\Omega_2(t) &= \frac{1}{2} \int_{0}^t \int_{0}^{t_1} \left[A(t_1),A(t_2)\right]\, dt_2dt_1 \ \ \ \mbox{and} \\
\Omega_3(t) &= \frac{1}{6} \int_{0}^t \int_{0}^{t_1} \int_{0}^{t_2} 
\left[A(t_1), \left[A(t_2),A(t_3)\right]\right] +  \left[A(t_3), \left[A(t_2),A(t_1)\right]\right]\,  dt_3dt_2dt_1.
\end{aligned}
\label{introduction:magnusterms}
\end{equation}
The straightforward evaluation of the $\Omega_j$ is nightmarishly expensive; however, a clever
technique which renders the calculations manageable is introduced in \cite{Iserles101}
and it paved the way for the development of a class of  numerical solvers 
which represent a fundamental matrix for 
(\ref{introduction:system}) over an interval $I$ via a collection of truncated
Magnus expansions.    While the entries of the $\Omega_j$ are slowly-varying whenever the entries of $A(t)$ are slowly-varying, the 
radius of convergence of the series  in (\ref{introduction:Magnus})
depends on the magnitude of the coefficient matrix $A(t)$, which is, in turn,
related to the magnitudes of the eigenvalues of $A(t)$.
Of course, this means that the number of Magnus expansions
which are needed to solve a given problem, and hence the cost of the method,
grows with the  magnitudes of the eigenvalues of $A(t)$.
See, for instance, \cite{Iserles102}, which gives for estimates
of the growth in the running time of Magnus expansion methods in 
the case of an equation of the form (\ref{introduction:two}) as 
a function of the magnitude of the coefficient $q$.

Nonetheless,  Magnus expansion methods 
are  much more efficient than standard solvers for ordinary
differential equations in the high-frequency regime.  Indeed, 
 exponential integrators which approximate Magnus expansions while avoiding the explicit calculation of 
commutators (those discussed in \cite{Blanes3}, for instance) appear to be the current state-of-the-art approach to solving
scalar ordinary differential equations of order three or higher. 
In our experiments, we compare our method against $4^{th}$ and $6^{th}$ order
``classical'' Magnus methods which explicitly make use of commutators, as well
as $4^{th}$ and $6^{th}$ order commutator-free quasi-Magnus exponential integrators.
Since the running time of our algorithm is largely independent
of frequency, our method is orders of magnitude faster  than Magnus-type methods in 
the high-frequency regime.  Perhaps surprisingly, we find that it is also faster 
even at quite low frequencies. We note, though, that Magnus
expansion methods are more general than our method in that they apply to
systems of linear ordinary differential equations and not just scalar equations.
Our experiment comparing our approach with Magnus-type methods
is described in  Subsection~\ref{section:experiments:2}.

We also compare our method with two specialized algorithms for second order equations:
the smooth deformation method of \cite{BremerPhase} (which was developed by one of the authors
of this paper) and  the ARDC method of \cite{agocs}.
These represent current state-of-the-art approaches to  solving  second order equations in the high-frequency regime.  
In the comparison made in Subsection~\ref{section:experiments:1},
we find that, despite its much greater generality, the algorithm of this paper is only slightly slower than
that of \cite{BremerPhase} and it is as much as 15 times faster than 
the ARDC method of \cite{agocs}.

The remainder of this article is organized as follows.  In Section~\ref{section:phase},
we discuss the results of \cite{BremerRokhlin} pertaining to the existence 
of slowly-varying phase functions for second order linear ordinary differential equations.
Section~\ref{section:levin} describes how the Levin principle can be exploited
to compute these slowly-varying phase functions.
In Section~\ref{section:algorithm}, we detail our numerical algorithm.
The results of numerical experiments demonstrating  the properties of our algorithm
are discussed in Section~\ref{section:experiments}.  These experiments include comparisons
with state-of-the-art methods for the special case of second order linear ordinary
differential equations and with Magnus-type exponential integrators.
We briefly comment on the algorithm of  this article and directions
for future work in Section~\ref{section:conclusion}.  Appendix~\ref{section:appendix}
details a standard adaptive spectral solver for ordinary differential equations which is used by our algorithm
and to construct reference solutions in our numerical experiments.

\end{section}

\begin{section}{Slowly-varying phase functions for second order equations}
\label{section:phase}

Here, we briefly discuss the results of  \cite{BremerRokhlin}, which pertain to
second order equations  of the form 
\begin{equation}
y''(t) + \omega^2 q_0(t) y(t) = 0,\ \ \ \ \ a < t < b,
\label{phase:two}
\end{equation}
with $q_0$ smooth and  positive.  Under these assumptions, the solutions of (\ref{phase:two}) 
are oscillatory, with the frequency of their oscillations controlled by the parameter 
$\omega$.   Analogous results hold when  $q_0$ is negative and the solutions
of (\ref{phase:two}) are combinations of rapidly increasing and decreasing
functions.  It is not obvious, however, how to apply the
argument of \cite{BremerRokhlin} to higher order scalar equations.
Nonetheless, there are strong indications, including relevant
 well-known results in asymptotic analysis (see, for instance, \cite{Wasov}) and experimental evidence,
that the situation for higher order scalar equations is similar.

If $y(t) = \exp(\psi(t))$ satisfies 
(\ref{phase:two}), then it can be trivially verified that  $\psi$ solves the Riccati equation
\begin{equation}
\psi''(t) + \left(\psi'(t)\right)^2 + \omega^2 q_0(t) = 0.
\label{phase:riccati}
\end{equation}
By inserting the expression $\psi(t) = i \alpha(t) + \beta(t)$  into (\ref{phase:riccati}),
we see that if   $\alpha$ and $\beta$  satisfy the system of equations
\begin{equation}
\left\{
\begin{aligned}
\beta''(t) +(\beta'(t))^2-(\alpha'(t))^2+\omega^2 q_0(t)&=0\\
 \alpha''(t) + 2\alpha'(t) \beta'(t) &= 0,
\end{aligned}
\right.
\label{phase:system}
\end{equation}
then $\psi$ solves (\ref{phase:riccati}).  The second equation in (\ref{phase:system})
admits the formal solution 
\begin{equation}
\beta(t) = -\frac{1}{2}\log(\alpha'(t)),
\label{phase:beta}
\end{equation}
so that $\psi$ can be written in the form
\begin{equation}
\psi(t) = i \alpha(t)-\frac{1}{2}\log\left(\alpha'(t)\right).
\end{equation}
Because of the close relationship 
between $\alpha$ and $\psi$, both are referred to as phase functions for (\ref{phase:two}).
Moreover, a bound  on the complexity of one readily gives a bound on the complexity of the other.

Inserting (\ref{phase:beta})  into the first equation in (\ref{phase:system}) yields 
\begin{equation}
\omega^2 q_0(t) - (\alpha'(t))^2 + \frac{3}{4} \left(\frac{\alpha''(t)}{\alpha'(t)}\right)^2
-\frac{1}{2}\frac{\alpha'''(t)}{\alpha'(t)}=0.
\label{phase:kummer}
\end{equation}
Equation (\ref{phase:kummer}) is known as  Kummer's equation, after E.~E.~Kummer, who studied it in \cite{Kummer}.
The theorem of  \cite{BremerRokhlin} applies when the function $p(x) = \widetilde{p}(t(x))$, where $\widetilde{p}(t)$ is defined via
\begin{equation}
\widetilde{p}(t) = 
\frac{1}{\omega^2 q_0(t)} \left(
\frac{5}{4}\left(\frac{q_0'(t)}{q_0(t)}\right)^2 - \frac{q_0''(t)}{q_0(t)}
\right)
\end{equation}
and $t(x)$ is the inverse function of 
\begin{equation}
x(t) = \int_a^t \sqrt{q_0(s)}\ ds,
\label{phase:xt}
\end{equation}
has a rapidly decaying Fourier transform.  More explicitly, the theorem
asserts that if the Fourier transform of $p$ satisfies a bound of the form
\begin{equation}
\left|\widehat{p}(\xi)\right|\leq \Gamma \exp\left(-\mu\left|\xi\right|\right),
\end{equation}
then there exist functions $\nu$ and $\delta$ such that 
\begin{equation}
\left|\nu(t)\right| \leq  \frac{\Gamma}{2\mu} \left(1 + \frac{4\Gamma}{\omega}\right) \exp(-\mu \omega),
\end{equation}
\begin{equation}
\left|\widehat{\delta}(\xi)\right| \leq \frac{\Gamma}{\omega^2}\left(1+\frac{2\Gamma}{\omega}\right)\exp(-\mu|\xi|)
\end{equation}
and
\begin{equation}
\alpha(t) = \omega \sqrt{q_0(t)} \int_a^t \exp\left(\frac{\delta(u)}{2}\right)\ du
\end{equation}
is a phase function for 
\begin{equation}
y''(t) + \omega^2 \left(q_0(t) + \frac{\nu(t)}{4\omega^2} \right) y(t) = 0.
\label{phase:1}
\end{equation}
Because the magnitude of $\nu$ decays exponentially fast in $\omega$,
 Equation~(\ref{phase:1}) is identical to (\ref{phase:two}) for the purposes of numerical
computation when $\omega$ is of even very modest size.  The definition of the function $p(x)$ is ostensibly quite complicated;
however, $p(x)$ is, in fact, simply a constant multiple of Schwarzian derivative of the inverse function
$t(x)$ of (\ref{phase:xt}).

This result ensures that for all values of $\omega$,
(\ref{phase:two}) admits a phase function which is slowly-varying.
In the low-frequency regime, when $\omega$ is small, it can be the case
that all phase functions for (\ref{phase:two}) oscillate, but they do so at
low frequencies because $\omega$ is small.  Once $\omega$ becomes sufficiently large,
the function $\nu$ is vanishingly small, and the phase
function associated with (\ref{phase:1}) is, at least for the purposes of numerical computation,
a slowly-varying phase function for the original equation (\ref{phase:two}).
Since $\nu$ decays exponentially fast in $\omega$, this happens at extremely modest frequencies.

Because of this phenomenon, in the low-frequency regime, 
the running time of numerical algorithms based on phase functions tend to grow with frequency.  However, once
a certain frequency threshold is reached, the complexity of the phase functions becomes essentially independent
of frequency, or even slowly decreasing with frequency.
This phenomenon can be clearly seen in all of the numerical experiments of this paper
presented in Section~\ref{section:experiments}.

\end{section}

\begin{section}{The Levin approach to solving nonlinear ordinary differential equations}
\label{section:levin}

In its original application to oscillatory integrals, Levin's principle was used
to construct slowly-varying solutions to inhomogeneous \emph{linear} ordinary differential
equations.  However, it can also be exploited to construct slowly-varying
solutions of  nonlinear ordinary differential equations, specifically the $(n-1)^{st}$ 
order Riccati equation.

When Newton's method is applied to the $(n-1)^{st}$ order Riccati equation corresponding
to (\ref{introduction:scalarode}), the result is a sequence of linearized equations of the form 
\begin{equation}
y^{(n-1)}(t) + p_{n-2}(t) y^{(n-2)}(t) + \cdots + p_1(t) y'(t) + p_0(t) y(t) = f(t).
\label{levin:inhom}
\end{equation}
Assuming the coefficients $q_0,\ldots,q_{n-1}$  
and the the initial guess used to initiate the Newton procedure are slowly-varying, the coefficients
$p_0,\ldots,p_{n-2}$ and the right-hand side $f$ appearing
in the first linearized equation of the form (\ref{levin:inhom})
which arises will also be slowly-varying.  According to the Levin principle
that equation admits slowly-varying solutions. If such a solution is used
to update the initial guess, then the second Newton iterate will also be slowly-varying
and the second linear inhomogeneous equation which arises will have slowly-varying coefficients
and a slowly-varying right-hand side.  Continuing in this fashion results in a series of linearized equations
of the form (\ref{introduction:inhom}), all of which have slowly-varying coefficients
and slowly-varying right-hand sides.  Consequently, a slowly-varying solution of the Riccati equation 
can be constructed via Newton's method as long as an appropriate slowly-varying initial guess is known.

Conveniently enough, there is an obvious mechanism for generating $n$ slowly-varying
initial guesses for the $(n-1)^{st}$ order Riccati equation.  In particular,
the eigenvalues $\lambda_1(t),\ldots,\lambda_n(t)$ of the matrix (\ref{introduction:scalarcoef}),
which are often used as low-accuracy approximations of solutions
of the Riccati equation in asymptotic methods, 
are suitable as  initial guesses for the Newton procedure.

Complicating matters slightly is the fact that the differential operator
\begin{equation}
D\left[y\right](t) = y^{(n-1)}(t) + p_{n-2}(t) y^{(n-2)}(t) + \cdots + p_1(t) y'(y) + p_0(t) y(t)
\label{levin:dop}
\end{equation}
appearing on the left-hand side of (\ref{levin:inhom}) admits a nontrivial nullspace which
can  contain rapidly-varying functions when one or more of the $p_j$ is of large magnitude.  
It is a  central observation of Levin-type methods, however, that 
 when (\ref{levin:inhom}) admits slowly-varying solutions
along with  rapidly-varying ones,  a slowly-varying solution can be accurately and rapidly computed
provided some case is taken.
In particular, as long as one uses a Chebyshev spectral collocation scheme which
is sufficient to resolve the coefficients $p_0,\ldots,p_{n-1}$ as well as the right-hand side $f$
and the resulting linear system is solved via a truncated singular value decomposition, a high-accuracy
approximation of a slowly-varying solution of (\ref{levin:inhom}) is obtained.
Critically, the discretization need not be sufficient to resolve
the  rapidly-varying solutions of (\ref{levin:inhom}) so that the cost of 
solving the equation depends only on the complexity of the desired slowly-varying
solution, rather than on the complexity of the rapidly-varying elements
of the nullspace of (\ref{levin:dop}).
Numerical evidence to this effect in the case $n=2$ is provided in 
\cite{LevinLi} and \cite{LiImproved}, and a detailed analysis is given 
 in \cite{SerkhBremerLevin}.

\end{section}

\begin{section}{Numerical Algorithm}
\label{section:algorithm}

In this section, we describe our method for the construction of 
a collection of slowly-varying phase functions $\psi_1,\ldots,\psi_n$
such that (\ref{introduction:phaserep}) is a basis in the space
of solutions of a nondegenerate equation of the form (\ref{introduction:scalarode})
with slowly-varying coefficients.   
Once these phase functions have been constructed, any reasonable
initial or boundary value problem for (\ref{introduction:scalarode})
can be easily solved.  Recall that we use $r_1,\ldots,r_n$ to denote the
first derivatives of the phase functions $\psi_1,\ldots,\psi_n$.

The algorithm operates in two stages, each of which is detailed in a
subsection below.  In the first stage, the Levin principle is used
to find the values  of $r_1,\ldots,r_n$ and their derivatives up to order $(n-2)$
at a point in the solution domain of the scalar equation.  
In the second stage, the Riccati equation corresponding
to (\ref{introduction:scalarode})  is solved  using these values 
as  initial conditions in order to calculate $r_1,\ldots,r_n$ and their derivatives
through order $(n-2)$ over the entire solution interval
and the phase functions  $\psi_1,\ldots,\psi_n$ are obtained by integrating $r_1,\ldots,r_n$.

Our algorithm takes as input the following:
\begin{enumerate}
\item
the interval $[a,b]$ over which the equation is given;

\item
an external subroutine for evaluating the coefficients $q_0,\ldots,q_{n-1}$ in (\ref{introduction:scalarode});

\item
a subinterval $[a_0,b_0]$ of $[a,b]$ over which the Levin procedure is to be applied and a point $\sigma$
in that interval;

\item
a point $\eta$ on the interval $[a,b]$ and the desired values $\psi_1(\eta),\ldots,\psi_n(\eta)$ for the phase
functions at that point;

\item
an integer $k$ which controls the order of the piecewise Chebyshev expansions used to represent the 
phase functions and their derivatives; and

\item
a parameter $\epsilon$  which specifies the desired accuracy 
for the solutions of the Riccati equation computed in the second stage of the algorithm.

\end{enumerate}

The output of our algorithm comprises $n^2$ piecewise Chebyshev expansions of order $(k-1)$,
representing the phase functions $\psi_1,\ldots,\psi_n$ and their derivatives
through order $(n-1)$.  To be entirely clear, by a  $(k-1)^{st}$ order piecewise Chebyshev 
expansions  on the interval $[a,b]$, we mean  a sum of the form
\begin{equation}
\begin{aligned}
&\sum_{i=1}^{m-1} \chi_{\left[x_{i-1},x_{i}\right)} (t) 
\sum_{j=0}^{k-1} \lambda_{ij}\ T_j\left(\frac{2}{x_{i}-x_{i-1}} t + \frac{x_{i}+x_{i-1}}{x_{i}-x_{i-1}}\right)\\
+
&\chi_{\left[x_{m-1},x_{m}\right]} (t) 
\sum_{j=0}^{k-1} \lambda_{mj}\ T_j\left(\frac{2}{x_{m}-x_{m-1}} t + \frac{x_{m}+x_{m-1}}{x_{m}-x_{m-1}}\right),
\end{aligned}
\label{global:chebpw}
\end{equation}
where $a = x_0 < x_1 < \cdots < x_m = b$ is a partition of $[a,b]$,
$\chi_I$ is the characteristic function on the interval $I$ and 
$T_j$ is the Chebyshev polynomial of degree $j$.
We note that the terms appearing in the first line of (\ref{global:chebpw}) involve
the characteristic function of a half-open interval, while that appearing
in the second involves  the characteristic function of a closed interval.
This ensures that exactly one term in  (\ref{global:chebpw}) is nonzero for each point $t$ in $[a,b]$.

\begin{subsection}{The Levin procedure}

In this first stage of the algorithm, the values of $r_1,\ldots,r_n$ and their derivatives
through order $(n-2)$ at the point $\sigma$ in the subinterval $[a_0,b_0]$ are calculated.
It proceeds as follows:

\begin{enumerate}

\item
Construct the $k$-point extremal Chebyshev grid $t_1,\ldots,t_{k}$ on the
interval $[a_0,b_0]$
and the corresponding $k \times k$ Chebyshev spectral differentiation matrix $D$.
The nodes are given by the formula
\begin{equation}
t_j = \frac{b_0-a_0}{2} \cos\left(\pi \frac{n-j}{n-1}\right) + \frac{b_0+a_0}{2}.
\label{oneint:nodes}
\end{equation}
The matrix $D$ takes the vector of values 
\begin{equation}
\left(
\begin{array}{c}
f(t_1)\\
f(t_2)\\
\vdots\\
f(t_k)\\
\end{array}
\right)
\end{equation}
of a Chebyshev expansion of the form
\begin{equation}
f(t) = \sum_{j=0}^{k-1} p_j T_j\left(\frac{2}{b_0-a_0} t + \frac{b_0+a_0}{b_0-a_0} \right)
\end{equation}
to the vector 
\begin{equation}
\left(
\begin{array}{c}
f'(t_1)\\
f'(t_2)\\
\vdots\\
f'(t_k)\\
\end{array}
\right)
\end{equation}
of the values of its derivatives at the nodes $t_1,\ldots,t_j$.

\item
Evaluate the coefficients $q_0,\ldots,q_{n-1}$ at the points $t_1,\ldots,t_k$
by calling the external subroutine supplied by the user.

\item
Calculate the values of $n$ initial guesses $r_{1},\ldots,r_{n}$ for the Newton procedure 
at the nodes $t_1,\ldots,t_k$ by first computing the eigenvalues  of the coefficient matrices
\begin{equation}
A_j = \left(
\begin{array}{cccccccc}
0 & 1 & 0 & \cdots & 0 & 0 \\
0 & 0 & 1 & \cdots & 0 & 0 \\
\vdots &  &  & \ddots &  & \vdots\\
0 & 0 & 0 & \cdots & 1 & 0 \\
0 & 0 & 0 & \cdots & 0 & 1 \\
-q_0(t_j) & -q_1(t_j) & -q_2(t_j) & \cdots&  -q_{n-2}(t_j) & -q_{n-1}(t_j)\\
\end{array}
\right)
\label{oneint:aj}
\end{equation}
for $j=1,\ldots,k$.
More explicitly, the eigenvalues of $A_j$ give the values of $r_1(t_j),\ldots,r_n(t_j)$.
The values of the first $(n-2)$ derivatives of $r_1,\ldots,r_n$ at the nodes $t_1,\ldots,t_k$
are then calculated through repeated application of the spectral differentiation matrix $D$.

\item
Perform Newton iterations in order to refine each of the initial guesses $r_1,\ldots,r_n$.
Because the general form of the Riccati equation is quite complicated, we illustrate
the procedure when $n=2$, in which case
the Riccati equation is 
\begin{equation}
r'(t) + (r(t))^2 + q_1(t) r(t) + q_0(t) = 0.
\end{equation}
In each iteration, we perform the following steps:

\begin{enumerate}

\item
Compute the residual 
\begin{equation}
\xi(t) = 
r'(t) + (r(t))^2 + q_1(t) r(t) + q_0(t) 
\end{equation}
of the current guess at the nodes $t_1,\ldots,t_k$.

\item
Form a spectral discretization of the linearized operator
\begin{equation}
L\left[\delta\right](t) = \delta'(t) + 2 r(t) \delta(t) + q_1(t) \delta(t).
\end{equation}
That is, form the $k\times k$ matrix
\begin{equation}
B = D + 
\left(\begin{array}{ccccc}
2r(t_1) + q_1(t_1) &                 &  &  \\
       & 2r(t_2) + q_1(t_2) &        &  & \\
       &        & \ddots &        &  \\  
       &         &         & 2r(t_k) + q_1(t_k)
\end{array}\right).
\end{equation}

\item
Solve the $k\times k$ linear system
\begin{equation}
B 
\left(
\begin{array}{c}
\delta(t_1)\\
\delta(t_2)\\
\vdots\\
\delta(t_k)
\end{array}
\right)
 = 
-\left(
\begin{array}{c}
\xi(t_1)\\
\xi(t_2)\\
\vdots\\
\xi(t_k)
\end{array}
\right)
\label{oneint:system}
\end{equation}
and update the current guess:
\begin{equation}
\left(
\begin{array}{c}
r(t_1)\\
r(t_2)\\
\vdots\\
r(t_k)
\end{array}
\right)
=
\left(
\begin{array}{c}
r(t_1)\\
r(t_2)\\
\vdots\\
r(t_k)
\end{array}
\right)
+
\left(
\begin{array}{c}
\delta(t_1)\\
\delta(t_2)\\
\vdots\\
\delta(t_k)
\end{array}
\right).
\end{equation}

\end{enumerate}

We perform a maximum of $8$ Newton iterations and the procedure is terminated if the value of
\begin{equation}
\max_{j=1,\ldots,k} \left|\delta(t_j)\right|
\end{equation}
is smaller than 
\begin{equation}
100 \epsilon_0\,, 
\max_{j=1,\ldots,k} \left|r(t_j)\right|,
\end{equation}
where $\epsilon_0\approx 2.220446049250313\times 10^{-16}$ denotes machine zero for the IEEE double precision number
system.

\item
We use Chebyshev interpolation to evaluate  $r_1,\ldots,r_n$, and their derivatives
of orders through $(n-2)$ at the point $\sigma \in [a_0,b_0]$.  These are the outputs
of this stage of the algorithm.
\end{enumerate}

Standard eigensolvers often produce inaccurate results in the case of 
matrices of the form (\ref{oneint:aj}), particularly when the entries
are of large magnitude.  Fortunately, there are specialized techniques
available for companion matrices, and the matrices appearing in (\ref{oneint:aj})
are simply the transposes of such matrices.  Our implementation of the procedure
of this subsection uses the backward stable and highly-accurate technique of \cite{AURENTZ1,AURENTZ2}
to compute the eigenvalues of the matrices (\ref{oneint:aj}).  

Care must also be taken when solving the linear system (\ref{oneint:system})
since the associated operator has a nontrivial nullspace.  Most of the time,
the discretization being used is insufficient to resolve any part of that nullspace,
with the consequence that the matrix $B$ is well-conditioned.  However,
when elements of the nullspace are sufficiently slowly-varying, they can be captured
by the discretization, in which case the matrix $B$ will have small singular values.
Fortunately, it is known that this does not cause numerical difficulties
in the solution of (\ref{oneint:system}), provided a truncated singular value
decomposition is used to invert the system.  Experimental evidence to this effect was presented in \cite{LevinLi,LiImproved}
and a careful analysis of the phenomenon appears in  \cite{SerkhBremerLevin}.
Because the truncated singular value decomposition is quite expensive, we actually use
a rank-revealing QR decomposition to solve the linear system (\ref{oneint:system})
in our implementation of the procedure of this subsection.  This was found to be about five times faster,
and it lead to no apparent loss in accuracy.

Rather than computing the eigenvalues of each of the matrices (\ref{oneint:aj}) in order
to construct initial guesses for the Newton procedure, one could accelerate
the algorithm slightly by computing the eigenvalues of only one  $A_j$ and use
the constant functions $r_1(t) = \lambda_1(t_j), \ldots, r(t) = \lambda_n(t_j)$ 
as initial guesses instead.  We did not make use of this optimization in our
implementation of the algorithm of this paper.
\end{subsection}

\begin{subsection}{Construction of the phase functions.}

Next, for each $j=1,\ldots,n$, the Riccati equation is solved using the value of $r_j(\sigma)$ to specify the desired solution.
These calculations are performed via the adaptive  spectral method described in Appendix~\ref{section:appendix}.
The parameters $k$ and $\epsilon$ are passed to that procedure.  Since most solutions of the Riccati equation are rapidly-varying
and we are seeking a slowly-varying solution, these problems are extremely stiff.  The solver of Appendix~\ref{section:appendix}
is well-adapted to such problems; however, essentially any solver for stiff ordinary differential equations would serve in its place.
The result is a collection of $n^2-n$ piecewise Chebyshev expansions of order $(k-1)$ representing the derivatives 
of the phase functions $\psi_1,\ldots,\psi_n$ of orders 1 through $(n-1)$.
Finally, spectral integration is used to construct $n$ additional piecewise Chebyshev expansions which
represent the phase functions $\psi_1,\ldots,\psi_n$ themselves.
 The particular  antiderivatives are determined by the values 
$\psi_1(\eta), \ldots, \psi_n(\eta)$ specified as inputs to the algorithm.

\end{subsection}

\end{section}

\begin{section}{Numerical experiments}
\label{section:experiments}

In this section, we present the results of numerical experiments which were conducted to illustrate the properties
of the method of this paper.  We implemented our algorithm in Fortran and compiled our code with  version 13.2.1 of the 
GNU Fortran compiler.    All experiments were performed on a single core of a  workstation computer equipped 
with an AMD 3995WX processor and 256GB of RAM.   No attempt  was made to parallelize our code.
The large amount of RAM was needed to calculate reference solutions
using a standard ODE solver.

Our algorithm calls for computing the eigenvalues of matrices of the
form (\ref{introduction:scalarcoef}).  Unfortunately, standard eigensolvers
lose significant accuracy when applied to many matrices of this type. However, 
because the  transpose of (\ref{introduction:scalarcoef}) is a companion matrix,
we were able to use the highly-accurate and backward stable algorithm of  \cite{AURENTZ1,AURENTZ2} 
for computing the eigenvalues of companion matrices to perform these calculations.

In all of our experiments,  the value of the parameter $k$, which determines the order of the 
Chebyshev expansions used to represent phase functions was taken to be $16$,
the particular antiderivatives $\psi_1,\ldots,\psi_n$ of the functions $r_1,\ldots,r_n$ 
were chosen through the requirement that  $\psi_1(0) = \psi_2(0) = \cdots = \psi_n(0) = 0$
and the Levin procedure was performed on the  subinterval $[0.0,0.1]$.  
The parameter $\epsilon$  which controls the accuracy of the obtained phase functions
was taken to be $10^{-12}$.

We tested the accuracy of the method of this paper by using it to calculate solutions to initial and boundary
value problems for scalar equations and comparing the results to reference solutions
constructed via the  standard adaptive spectral method described in Appendix~\ref{section:appendix}.
Because the condition numbers of these initial and boundary value problems for (\ref{introduction:scalarode}) 
grow with frequency, the accuracy of any numerical
method used to solve them is expected to deteriorate with increasing  frequency.    In the case of our algorithm,
the phase functions themselves are calculated to high precision, but their magnitudes
increase with frequency and accuracy is lost when the phase functions  are exponentiated.
One implication is of this is that calculations which involve only the phase functions
and not the solutions of the scalar equation can be performed to high accuracy.
The article \cite{BremerZeros}, for example,  describes
a scheme of this   type for rapidly computing the zeros of solutions of second order linear ordinary
differential equations to extremely high accuracy.

To account for the vagaries of modern computing environments, all reported times were obtained 
by averaging the cost of each calculation over either 1,000 runs.


\begin{subsection}{Comparison with two specialized methods for second order equations}
\label{section:experiments:1}

We first compared the performance of the Levin-type method of this paper
with the smooth deformation scheme of \cite{BremerPhase} developed by one of this paper's  authors,
and with the ARDC method of \cite{agocs}.

For each $\nu=2^0,2^1,2^2,\ldots,2^{20}$ and each of the three methods considered,
we solved Legendre's differential equation
\begin{equation}
(1-t^2) y''(t) - 2t y'(t) + \nu(\nu+1) y(t) = 0
\end{equation}
in order to obtain the Legendre polynomial $P_\nu$ of degree $\nu$ over the interval $[0.0,0.999]$.
The algorithm  of \cite{agocs} makes it somewhat difficult
to evaluate solutions at arbitrary points inside the solution domain, so we
settled for measuring the error in each obtained solution by comparing its
value at $t=0.999$ against the known value of $P_\nu(0.999)$. 

We used the implementation of the method of \cite{BremerPhase} available at:
\begin{center}
\url{https://github.com/JamesCBremerJr/Phase-functions}
\end{center}
We used an implementation of the ARDC method designed specifically for solving Legendre's
differential equation which was suggested to us by one of the authors of \cite{agocs}. 
It is available at:
\begin{center}
https://github.com/fruzsinaagocs/riccati/tree/legendre-improvements
\end{center}
The more general implementation of the ARDC method used in the experiments of 
 \cite{agocs}, which does not perform as well in this experiment,
can be found at:
\begin{center}
\url{https://github.com/fruzsinaagocs/riccati}
\end{center}

The input parameters for the algorithms of \cite{BremerPhase} and \cite{agocs} were set as follows.
For the method of \cite{BremerPhase},  we set the parameter $k$ controlling the order of the piecewise
Chebyshev expansions used to represent phase functions to be $16$, and took the parameter $\epsilon$
specifying the desired accuracy for the phase functions to be $10^{-12}$.
For \cite{agocs}, we used the default parameters provided by the authors' code.

\begin{figure}[t!!!!!!!!!!!!!!!!!!!!]
\hfil
\includegraphics[width=.32\textwidth]{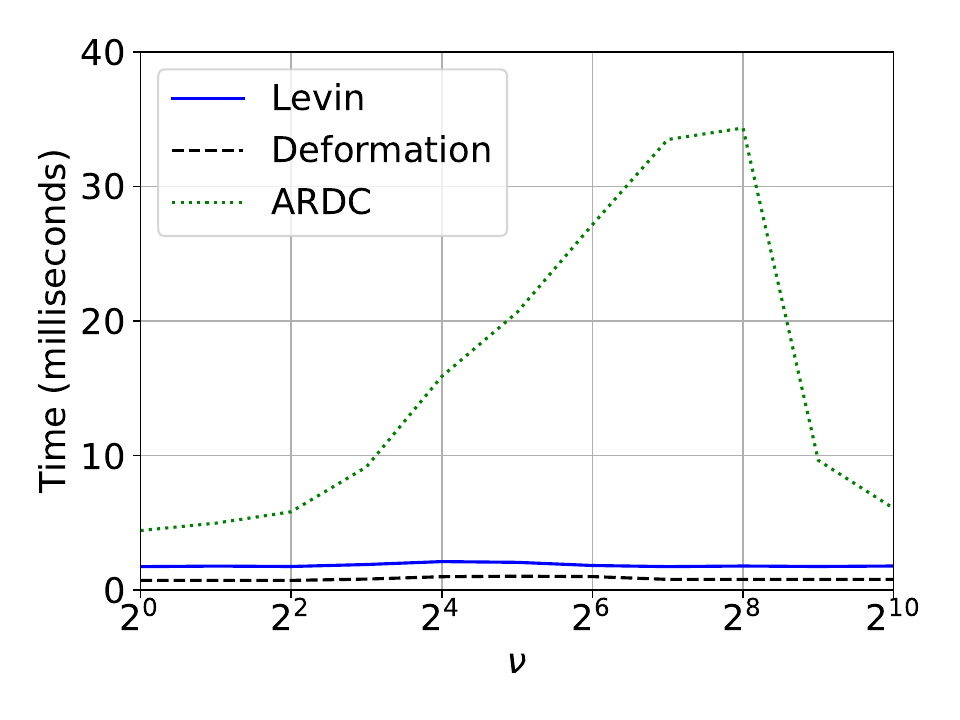}
\hfil
\includegraphics[width=.32\textwidth]{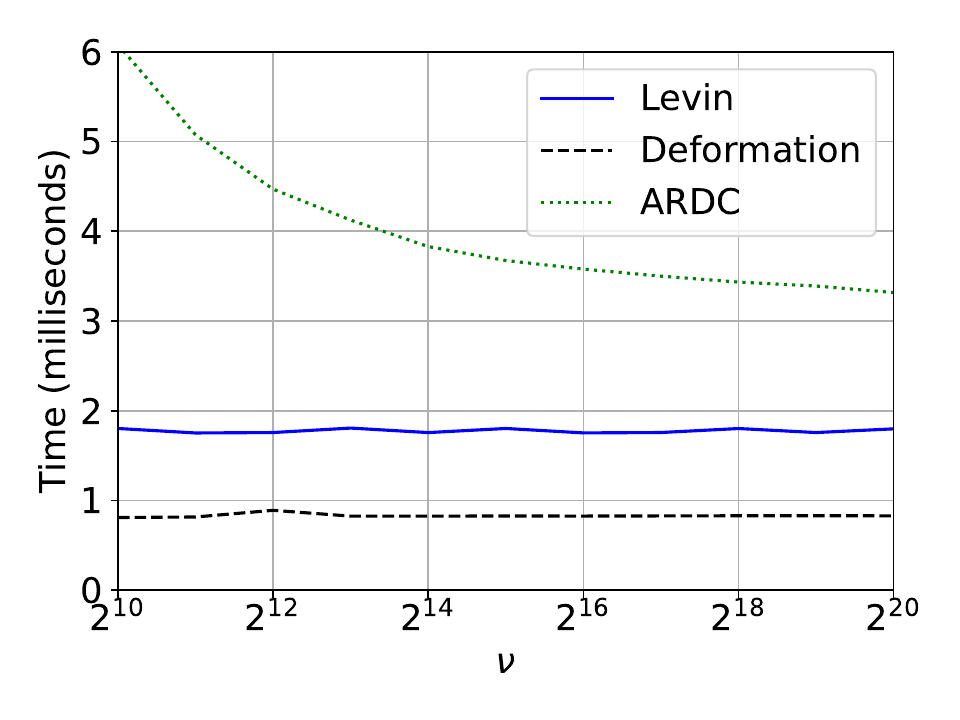}
\hfil
\includegraphics[width=.32\textwidth]{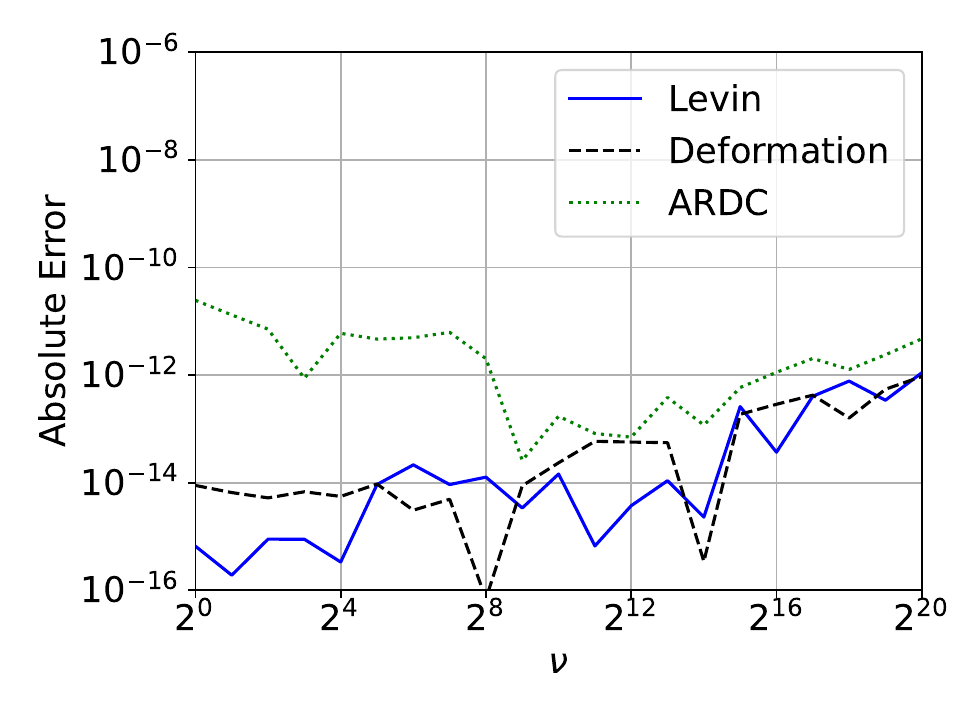}
\hfil

\caption{The results of the experiment of Subsection~\ref{section:experiments:1}
in which the Levin-type method of this paper, the smooth deformation scheme of \cite{BremerPhase}
and  the ARDC method of \cite{agocs} are compared.  The left-most plot gives the time
required by each algorithm as a function of $\nu$, but only for the low-frequency regime.
The middle plot gives the time required by each algorithm in the high-frequency regime.
The plot on the right shows the  absolute error in the value of the Legendre $P_\nu(0.999)$ obtained by each algorithm
as a function of $\nu$.   
}
\label{experiments1:figure1}
\end{figure}

Figure~\ref{experiments1:figure1} presents the results of this experiment. 
We observe that the method of this paper achieves similar accuracy to that of \cite{BremerPhase},
but is a bit slower.
Although \cite{agocs} claims that ARDC achieves
a ten times speed improvement over the method of \cite{BremerPhase}, we have not found this to be the case.
At frequencies below $2^9$, the ARDC method is both noticeably slower and less accurate than
both the other methods.  For example, when $\nu=2^8$, 
the algorithm of this paper takes around 1.8 milliseconds and achieves 13 digits of accuracy,  
that of \cite{BremerPhase} takes approximately 0.81 milliseconds and achieves 15 digits of accuracy
while  the ARDC method takes more than 30 milliseconds and obtains only 11 digits of accuracy. 
In particular, ARDC can be as much as 15 times slower than the method of this paper and 30 times
slower than the algorithm of \cite{BremerPhase}.
At higher frequencies, ARDC achieves similar levels of accuracy to \cite{BremerPhase} and the method of this paper, 
but it is more than a factor of two slower than the algorithm of this paper and more than
a factor of three slower than the method of \cite{BremerPhase}.
The discrepancy between results reported in \cite{agocs} and the results of this experiment
appears to be attributable to the use of an unoptimized, highly inefficient implementation of \cite{BremerPhase}
in the experiments of \cite{agocs}.

As explained in Section~\ref{section:phase}, in the low-frequency regime,  the running times of all three methods increase with $\omega$.
However, once a certain frequency threshold is reached, the running times
decrease rapidly and then become essentially independent of frequency, or even continue to decrease
slowly as functions of $\omega$.  
We note that, in our plots, this phenomenon is more apparent in the case of the ARDC  method because of the 
much greater cost of that algorithm in the low-frequency regime.

\end{subsection}

%
%
%

\begin{subsection}{Comparison with Magnus-type exponential integrators}
\label{section:experiments:2}

In our second experiment, we compared the performance of our algorithm with that of 
four methods based on Magnus-type exponential integrators.  We use MG4 to refer to the 
 $4^{th}$ order Magnus exponential integrator given by (2.9) in \cite{Iserles103};
MG6 denotes the $6^{th}$ order Magnus exponential integrator specified by (3.10) in
\cite{Blanes1};
we use CF4 to refer to  $4^{th}$ order two exponential commutator-free quasi-Magnus exponential integrator listed
in Table~2 of \cite{Blanes3};
and CF6 is the first of the $6^{th}$ order five exponential commutator-free quasi-Magnus exponential integrators
listed in  Table~3 of \cite{Blanes3}.  
\begin{figure}[b!]
\begin{center}
\includegraphics[width=.50\textwidth]{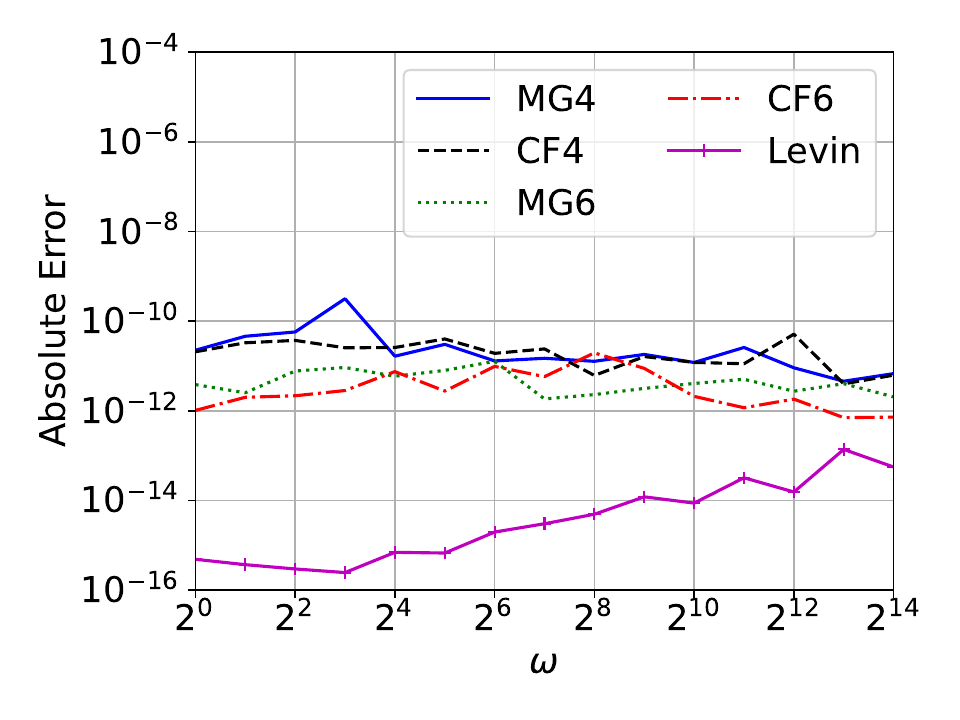}
\end{center}
\caption{The errors in the solutions of the initial value problem of Subsection~\ref{section:experiments:2}
obtained via four Magnus-type exponential integrator methods and the Levin-type algorithm of this paper. 
}
\label{experiments2:figure1}
\end{figure}

The performance of exponential integrator methods depends critically on proper step length control.  
In order to give every possible benefit to the methods we compare our scheme against, we 
use the following two-phased approach.  In the first phase, which was not timed, we determined a sequence
of appropriate step sizes via a greedy algorithm.  More explicitly, at each step, we started
with a  large step size $h$ and repeatedly reduced it by a factor of $0.95$ until 
an estimate of the local error fell bellow $\epsilon = 10^{-12}$.  The local error estimate
was obtained by taking two steps of length $h/2$ in order to produce a (hopefully)
superior approximation of the value of the solution at the terminal point.
In the second phase, the equation was solved using the precomputed sequence of step lengths.
It is only the second phase of the calculation which was timed.

For each $\omega=2^0,2^1,\ldots,2^{14}$ and each of the five methods, we  solved the differential equation
\begin{equation}
y'''(t) + q_2(t) y''(t) + q_1(t) y'(t) + q_0(t) y(t) = 0,
\label{experiments2:1}
\end{equation}
where
\begin{equation}
\begin{aligned}
q_0(t) &=-\frac{\omega  \left(e^t \omega -i\right) (\cos (8 t)+3) \left(\left(t^2+1\right) \cos (3 t)-i \omega \right)}{t^2+1}\\
q_1(t) &=\frac{\omega  \left(-\left(\omega +i \left(t^2+1\right)\right) \cos (8 t)+e^t \omega  \left(3 t^2+\left(t^2+1\right) \cos (8 t)+4\right)-3 i t^2-3 \omega -4 i\right)}{t^2+1}+\\
&\ \ \ \ \cos (3 t) \left(i \left(e^t-3\right) \omega -i \omega  \cos (8 t)+1\right)\ \ \ \mbox{and}
\\
q_2(t) &=i \left(\frac{1}{t^2+1}-e^t+3\right) \omega +i \omega  \cos (8 t)-\cos (3 t)-1,
\end{aligned}
\label{experiments2:2}
\end{equation}
%
%
over the interval $[0,0.1]$ subject to the conditions
\begin{equation}
y(0) = 1,\ \ \ y'(0) = i \omega \ \ \ \mbox{and}\ \ \ y''(0) = (i\omega)^2.
\label{experiments2:3}
\end{equation}
The eigenvalues of the coefficient matrix corresponding to Equation~(\ref{experiments2:1}) are
\begin{equation}
\lambda_1(t) = 1+i e^t \omega,\ \ \ 
\lambda_2(t) = \cos (3 t)-\frac{i \omega }{t^2+1}\ \ \ \mbox{and} \ \ \ 
\lambda_3(t) = -i \omega  (\cos (8 t)+3).
\label{experiments2:4}
\end{equation}
As in the case of the experiment of the last section, owing to the difficulty
of computing solutions at arbitrary points using step methods,
we assessed the accuracy of the obtained solutions by measuring
the absolute error in their values at the endpoint $t=0.1$ of the solution domain only.
Moreover, we only considered values of $\omega$ up to $2^{14}$ because the cost 
of finding  appropriate step sizes becomes excessive for larger  values of $\omega$.

\begin{table}[b!!]
\begin{center}
\begin{tabular}{cccccc}
\toprule
$\omega$ & MG4 & CF4 & MG6 & CF6 & Levin\\
\midrule
$2^{  0}$ & 2.79\e{-03} & 3.77\e{-03} & 9.70\e{-04} & 9.89\e{-04} & 6.88\e{-04}  \\
$2^{  1}$ & 3.72\e{-03} & 4.96\e{-03} & 1.46\e{-03} & 1.46\e{-03} & 7.07\e{-04}  \\
$2^{  2}$ & 7.42\e{-03} & 8.97\e{-03} & 2.95\e{-03} & 2.45\e{-03} & 7.35\e{-04}  \\
$2^{  3}$ & 1.51\e{-02} & 1.47\e{-02} & 5.42\e{-03} & 3.44\e{-03} & 8.91\e{-04}  \\
$2^{  4}$ & 2.55\e{-02} & 2.44\e{-02} & 9.71\e{-03} & 6.42\e{-03} & 7.57\e{-04}  \\
$2^{  5}$ & 4.42\e{-02} & 4.59\e{-02} & 1.94\e{-02} & 1.23\e{-02} & 7.60\e{-04}  \\
$2^{  6}$ & 7.79\e{-02} & 8.07\e{-02} & 3.48\e{-02} & 2.13\e{-02} & 7.61\e{-04}  \\
$2^{  7}$ & 1.35\e{-01} & 1.40\e{-01} & 6.46\e{-02} & 3.99\e{-02} & 7.62\e{-04}  \\
$2^{  8}$ & 2.48\e{-01} & 2.48\e{-01} & 1.13\e{-01} & 7.21\e{-02} & 7.63\e{-04}  \\
$2^{  9}$ & 4.35\e{-01} & 4.40\e{-01} & 2.14\e{-01} & 1.31\e{-01} & 7.43\e{-04}  \\
$2^{ 10}$ & 7.61\e{-01} & 7.59\e{-01} & 3.95\e{-01} & 2.41\e{-01} & 7.42\e{-04}  \\
$2^{ 11}$ & 1.35\e{+00}  & 1.30\e{+00}  & 6.96\e{-01} & 4.27\e{-01} & 7.41\e{-04}  \\
$2^{ 12}$ & 2.25\e{+00}  & 2.25\e{+00}  & 1.29\e{+00}  & 7.93\e{-01} & 7.42\e{-04}  \\
$2^{ 13}$ & 3.88\e{+00}  & 3.95\e{+00}  & 2.27\e{+00}  & 1.41\e{+00}  & 7.40\e{-04}  \\
$2^{ 14}$ & 6.78\e{+00}  & 7.02\e{+00}  & 4.36\e{+00}  & 2.74\e{+00}  & 7.41\e{-04}  \\
\bottomrule\\
\end{tabular}

\end{center}
\caption{The time, in second, required by four Magnus-type exponential integrator methods and the Levin-type algorithm of this paper
to solve the initial value problem of Subsection~\ref{section:experiments:2}.}
\label{experiments2:table1}
\end{table}

Figure~\ref{experiments2:figure1} and Table~\ref{experiments2:table1} give the results.  
We observe that all of the methods achieve reasonably accuracy given the requested
level of precision.  Not surprisingly, given the difference in the asymptotic behaviour of the running time
of these algorithms with respect to frequency,
the algorithm of this paper  is orders of magnitude faster than
the exponential integrator methods at high frequencies.  In fact,
when $\omega=2^{14}$, our approach is more than $3,000$ times faster than the 
most efficient of the exponential integrator methods.  What is perhaps surprising, is that
the algorithm of this paper is faster than the various exponential integrator methods
even at very low frequencies.  This is indicative of the fact that, even in the low-frequency
regime, phase functions are not much more expensive to represent than the solutions of the 
scalar equation itself.

\end{subsection}

\begin{figure}[h!]
\hfil
\includegraphics[width=.32\textwidth]{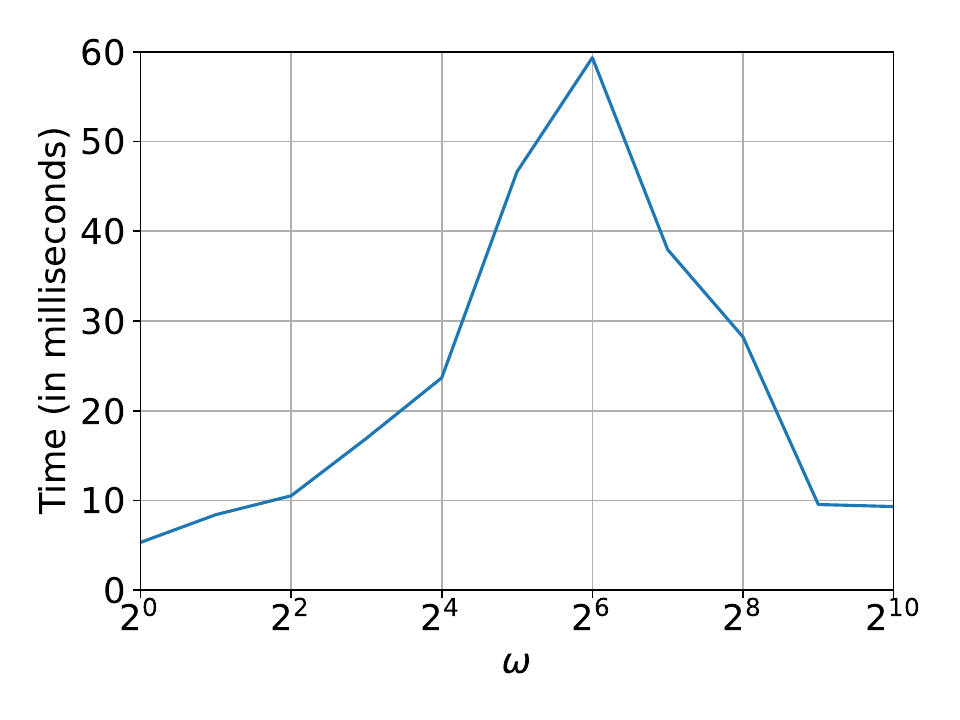}
\hfil
\includegraphics[width=.32\textwidth]{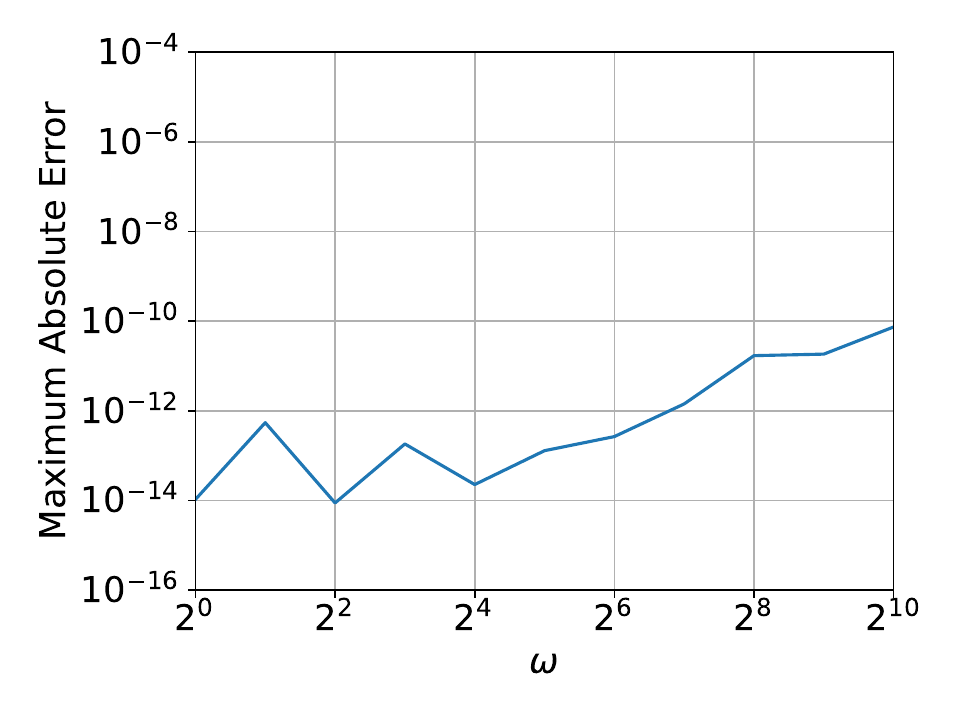}
\hfil
\includegraphics[width=.32\textwidth]{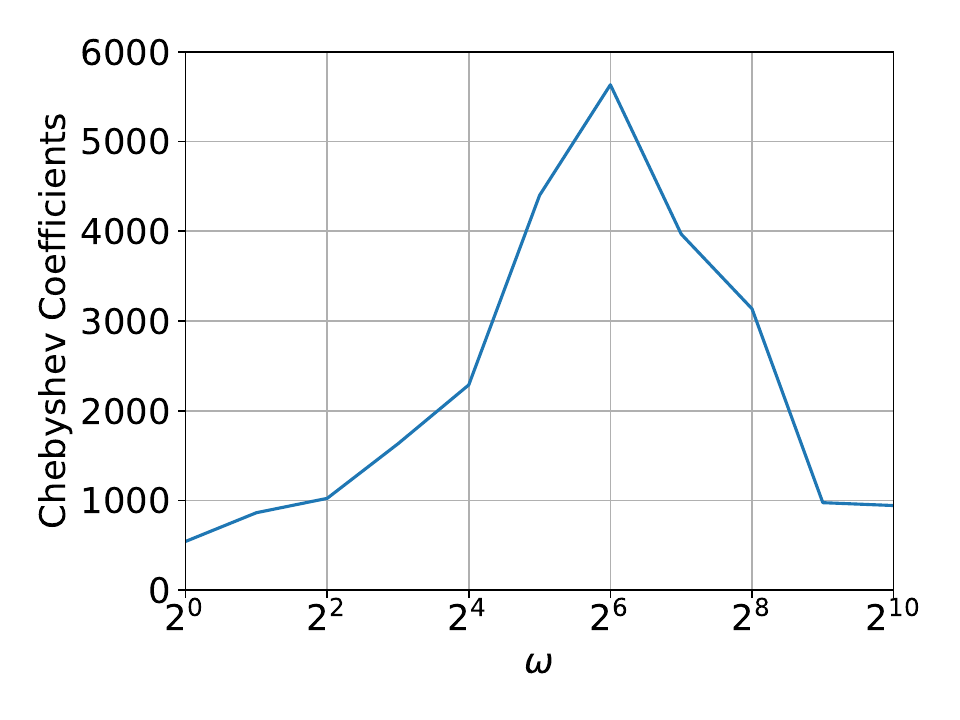}
\hfil

\hfil
\includegraphics[width=.32\textwidth]{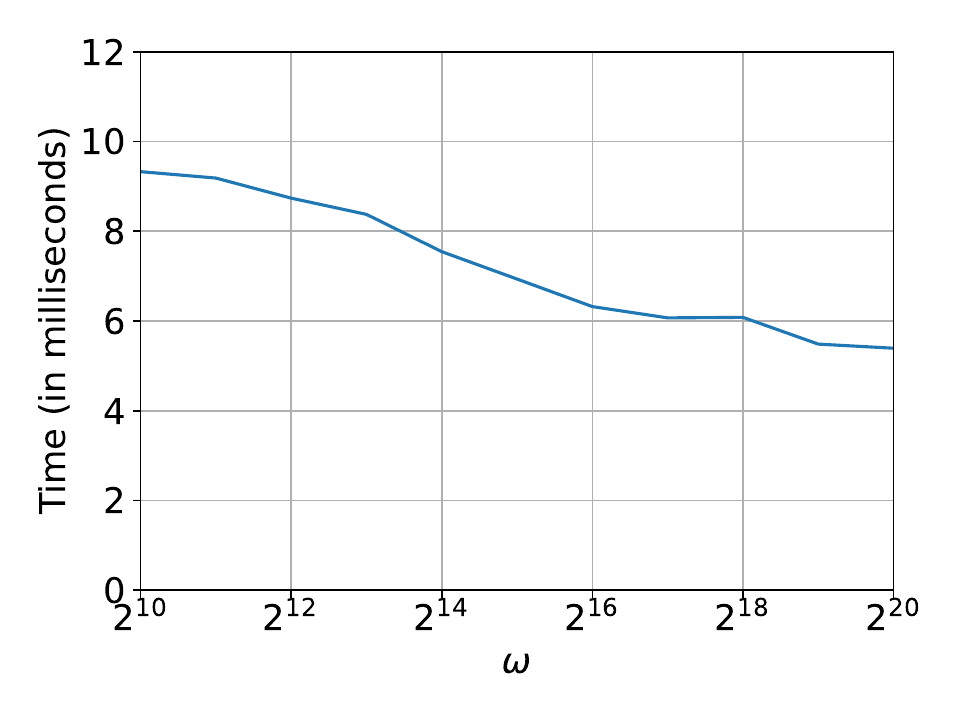}
\hfil
\includegraphics[width=.32\textwidth]{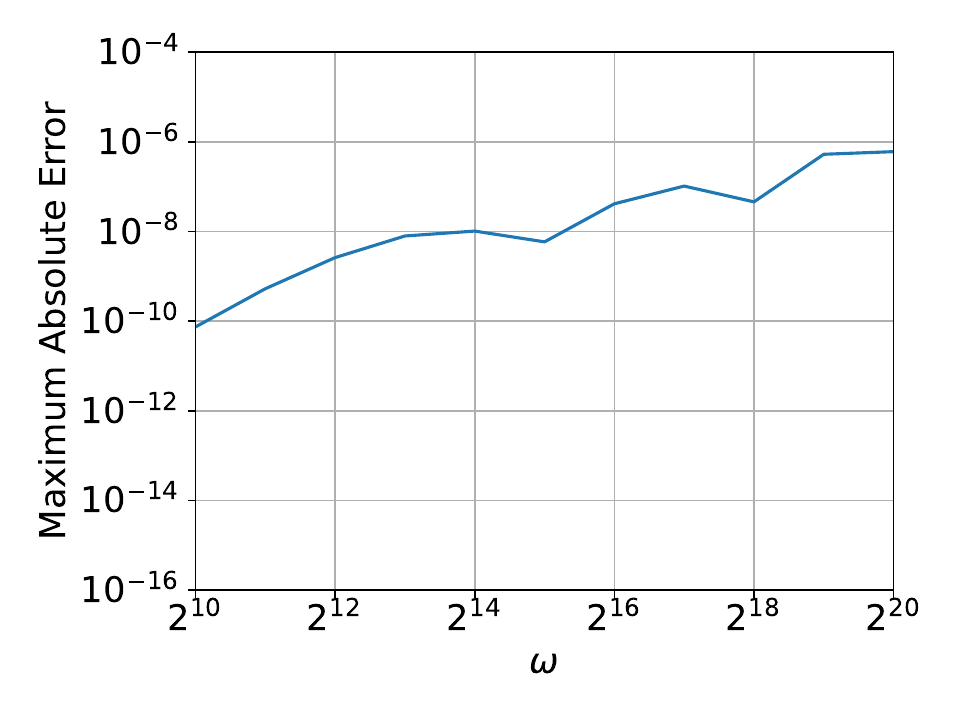}
\hfil
\includegraphics[width=.32\textwidth]{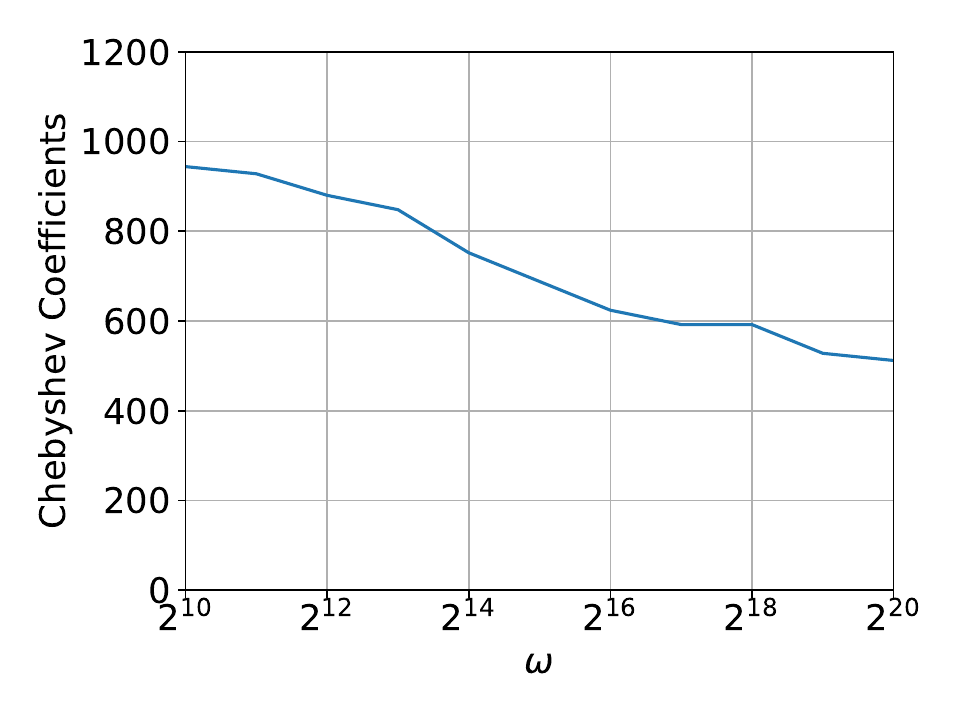}
\hfil

\caption{The results of the experiments of Subsection~\ref{section:experiments:3}.
The plot at top left gives the running time of the method of this paper
in the low-frequency regime.
The top-middle plot gives reports the absolute error in the solution of the boundary value problem for (\ref{experiments3:1})
in the low-frequency regime.
The plot at top right shows  the total number  of piecewise Chebyshev coefficients required to represent the slowly-varying phase functions,
 again in the low-frequency regime.
The plots on the bottom row provide the same information, but in the high-frequency regime.
}
\label{experiments3:figure1}
\end{figure}

\begin{figure}[h!!!!!!!]
\hfil
\includegraphics[width=.28\textwidth]{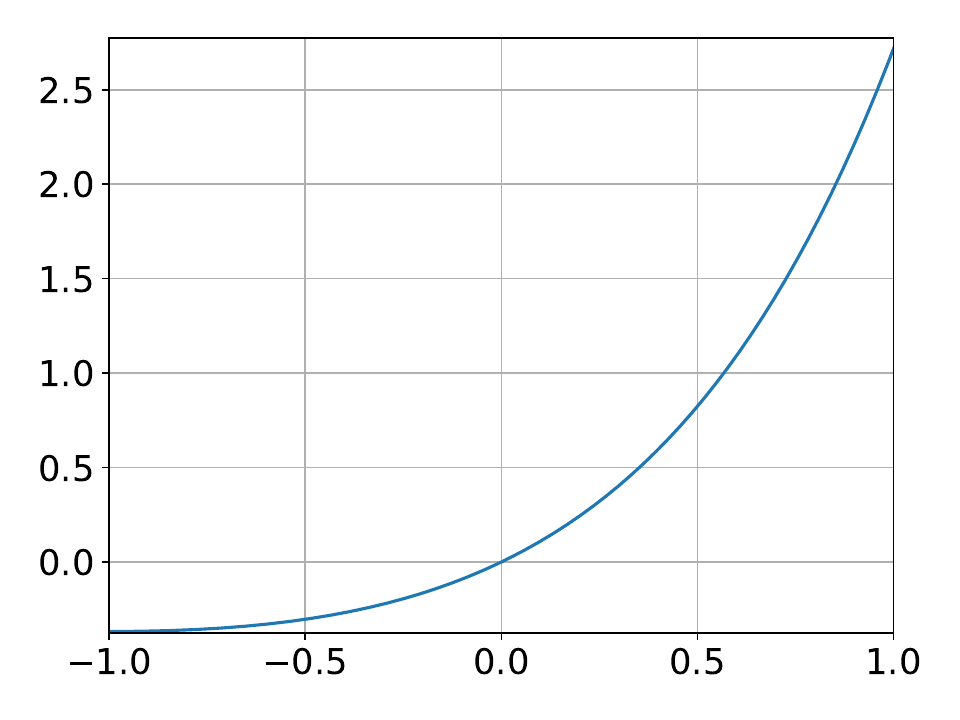}
\hfil
\includegraphics[width=.28\textwidth]{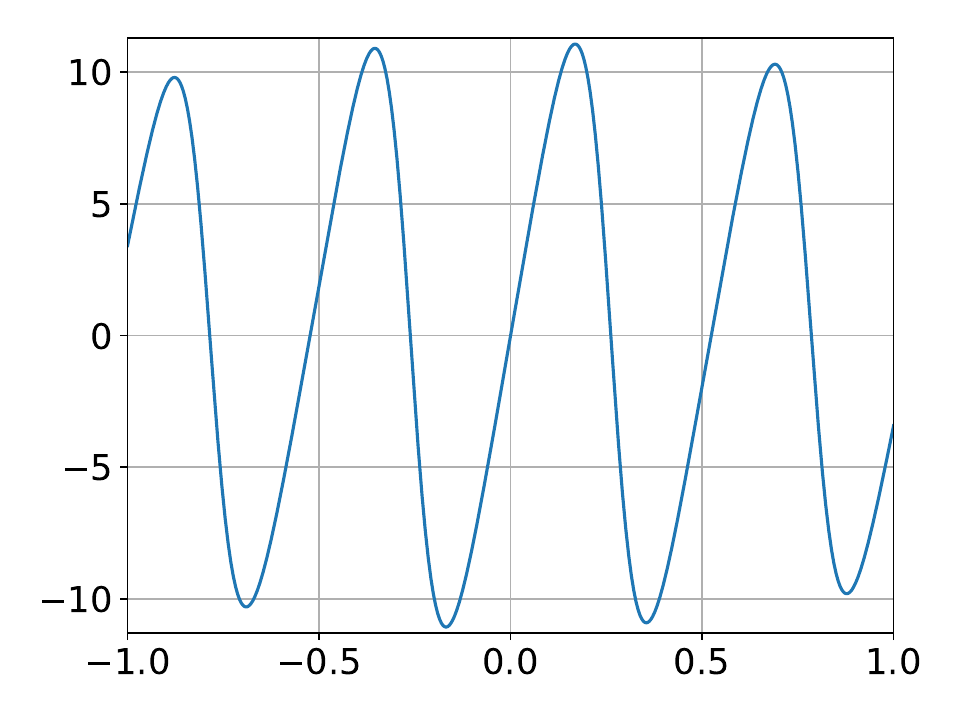}
\hfil
\includegraphics[width=.28\textwidth]{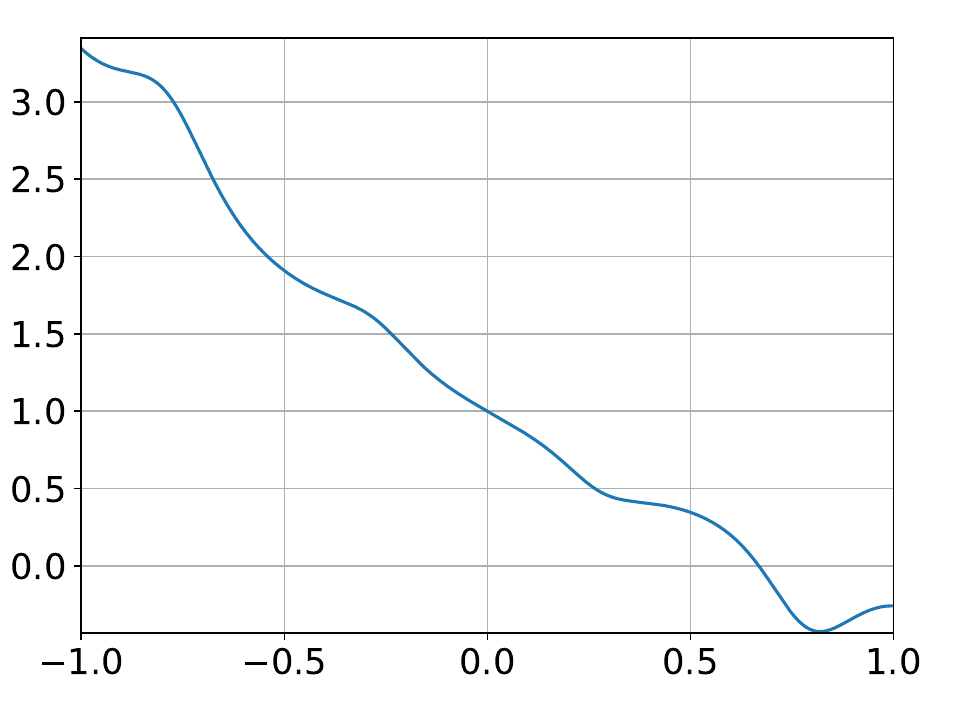}
\hfil

\hfil
\includegraphics[width=.28\textwidth]{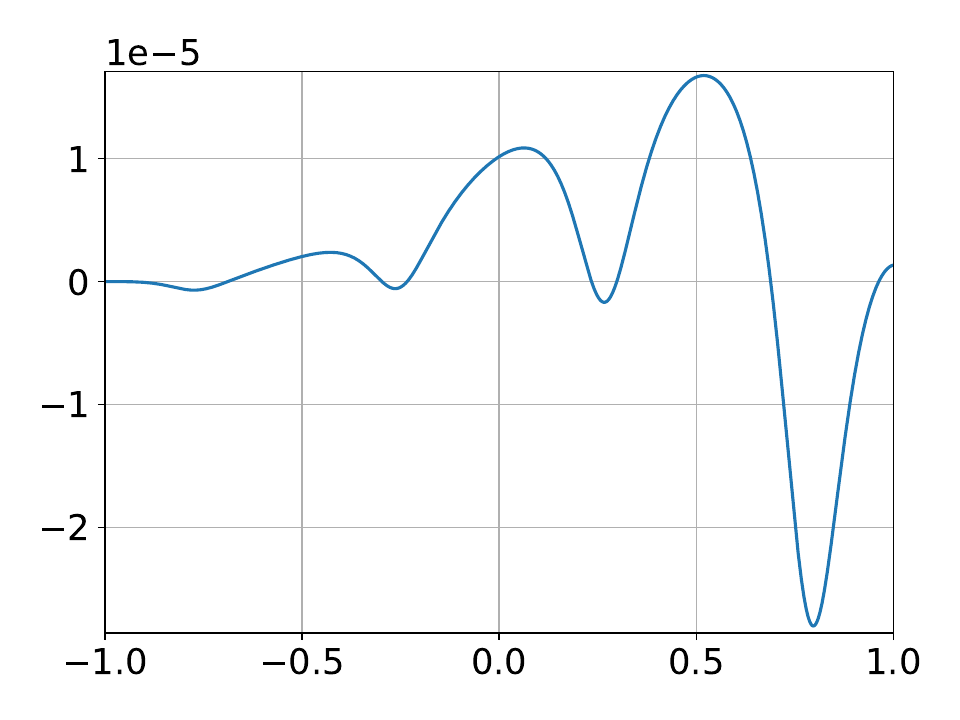}
\hfil
\includegraphics[width=.28\textwidth]{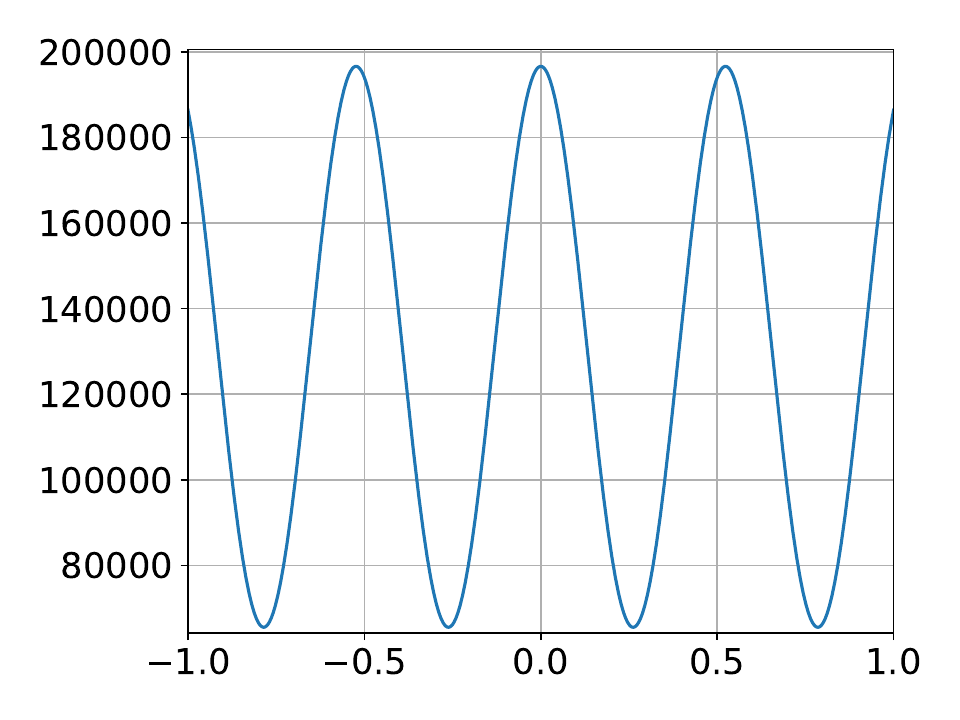}
\hfil
\includegraphics[width=.28\textwidth]{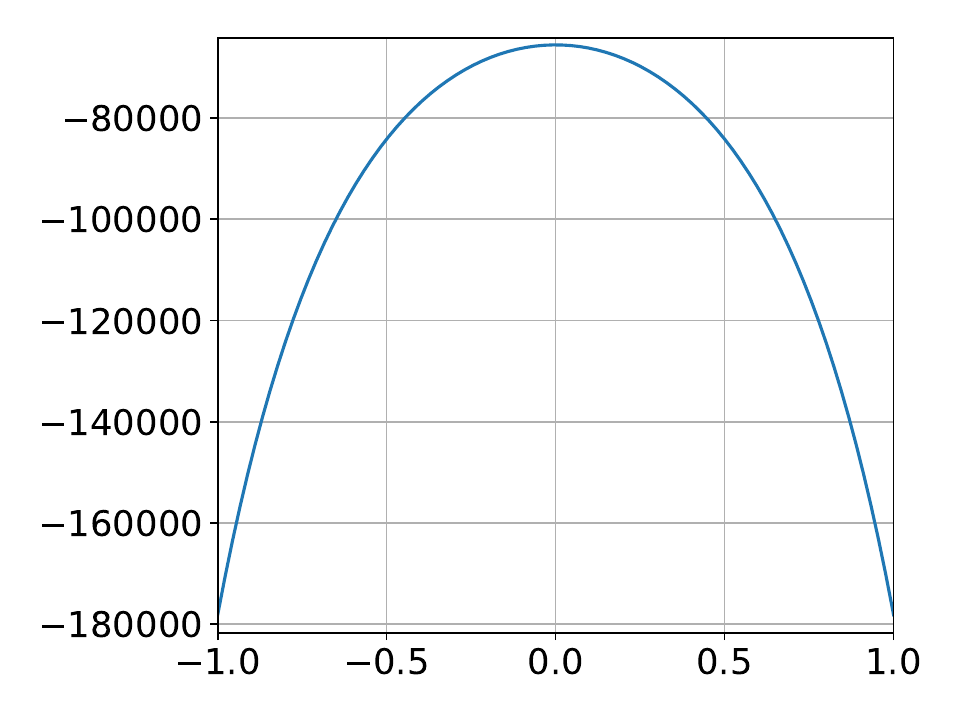}
\hfil

\caption{The derivatives of the three slowly-varying phase functions produced
by applying the method of this paper to Equation~(\ref{experiments3:1}) of Subsection~\ref{section:experiments:3}
when the parameter $\omega$ is equal to $2^{16}$.  Each column corresponds to one of the phase functions,
with the real part appearing in the first row and the imaginary part in the second.}
\label{experiments3:figure2}
\end{figure}

%
%
\begin{subsection}{A boundary value problem for a third order equation}
\label{section:experiments:3}

In the experiment described in this section, we considered the equation
\begin{equation}
y'''(t) + q_2(t) y''(t) + q_1(t) y'(t) + q_0(t) y(t) = 0,
\label{experiments3:1}
\end{equation}
where
\begin{equation}
\begin{aligned}
q_0(t) &=-i e^t t \omega  \left(e^t-i e^{t^2} \omega \right) (\cos (12 t)+2),\\
q_1(t) &=e^{t^2} \omega  \left(2 \omega -i e^t t\right)+\omega  \left(e^{t^2} \omega +i e^t (t+1)\right) \cos (12 t)+e^t \left(e^t t+2 i (t+1) \omega \right)
\ \ \ \mbox{and}
\\
q_2(t) &=i e^{t^2} \omega -i \omega  \cos (12 t)-e^t (t+1)-2 i \omega.
\end{aligned}
\label{experiments3:2}
\end{equation}
The eigenvalues of the coefficient matrix corresponding to (\ref{experiments3:1}) are 
\begin{equation}
\begin{aligned}
\lambda_1(t) &=i \omega  (\cos (12 t)+2), \ \ \ 
\lambda_2(t) &=t e^t  \ \ \ \mbox{and} \ \ \ 
\lambda_3(t) &=e^t-i e^{t^2} \omega.
\end{aligned}
\end{equation}
For each $\omega=2^0,2^1,\ldots,2^{20}$, we used our algorithm to solve (\ref{experiments3:1}) over the interval
$[-1,1]$ subject to the conditions
\begin{equation}
y(-1) =  y(1) = 1 \ \ \ \mbox{and} \ \ \ y'(-1) = 0.
\label{experiments3:3}
\end{equation}
We measured the absolute error in each resulting solution at 10,000 equispaced points
in the interval $[-1,1]$ via comparison with a reference solution constructed
using the solver of Appendix~\ref{section:appendix}.  
The results are given in Figure~\ref{experiments3:figure1} while
Figure~\ref{experiments3:figure2} contains plots of 
the derivatives of the three slowly-varying phase functions produced
by applying the method of this paper to Equation~(\ref{experiments3:1})
when $\omega=2^{16}$.
As expected, the running time of the method of this paper increases 
until a certain frequency threshold is passed, at which point it
falls precipitously before becoming slowing decreasing.  
The maximum observed absolute error in the solution grows consistently with $\omega$, which is
as expected considering that the condition number of the problem deteriorates with increasing
frequency.
For all values of $\omega$ greater than or equal to $2^{9}$, less than 10 milliseconds
was required to solve the boundary value problem  and fewer than 1,000 Chebyshev coefficients were
needed to represent the phase functions.  No more than 60 milliseconds and 6,000 coefficients
were required in the worst case.  
The frequency $\Omega$ of the problems considered increased from approximately 3.9 when $\omega=1$ to
roughly 4,100,531 when $\omega=2^{20}$.

\end{subsection}

%
%
\begin{subsection}{An initial value problem for a fourth order equation}
\label{section:experiments:4}

In this experiment, we considered the linear scalar ordinary differential equation
\begin{equation}
y''''(t) + q_3(t) y'''(t) + q_2(t) y''(t) + q_1(t) y'(t) + q_0(t) y(t) = 0
\label{experiments4:1}
\end{equation}
whose coefficient matrix has eigenvalues
\begin{equation}
\begin{aligned}
\lambda_1(t) =\frac{t}{2}+i e^{t^2} \omega,\ \ \
\lambda_2(t) =\frac{i \omega }{t^2+2}+e^{i t},\ \ \
\lambda_3(t) =\cos (3 t)\ \ \ \mbox{and}\ \ \
\lambda_4(t) =-i \left(t^2+1\right) \omega.
\end{aligned}
\label{experiments4:2}
\end{equation}
Formulas for the coefficients $q_0$, $q_1$, $q_2$ and $q_3$  are too unwieldy
to reproduce here, but they can be easily calculated from (\ref{experiments4:2})
using a computer algebra system.
For each $\omega=2^0,2^1,\ldots,2^{20}$, we used the algorithm of this paper to solve (\ref{experiments4:1}) over the interval
$[-1,1]$ subject to the conditions
\begin{equation}
y(0) = 1, \ \ \ y'(0) = i \omega, \ \ \ y''(0) = (i \omega)^2 \ \ \ \mbox{and} \ \ \   y'''(0) = (i \omega)^3.
\label{experiments4:3}
\end{equation}
We measured the absolute error in each resulting solution at 10,000 equispaced points
in the interval $[-1,1]$ via comparison with a reference solution constructed
using the  solver of Appendix~\ref{section:appendix}.  
The results are given in Figure~\ref{experiments4:figure1}.
We observe that for all $\omega$ greater than or equal to $2^9$, 
fewer than $8$ milliseconds was required to solve the problem and less
than $250$ piecewise Chebyshev coefficients were required to represent the phase functions.
In the worst case, when $\omega=2^6$, the solver took around 92 milliseconds and 3,200 piecewise Chebyshev coefficients
were needed.  The frequency $\Omega$ of the problems considered ranged from around 2.97 when $\omega=1$
to approximately $3,067,403$ when $\omega=2^{20}$.

\begin{figure}[h!!!!!!!]
\hfil
\includegraphics[width=.32\textwidth]{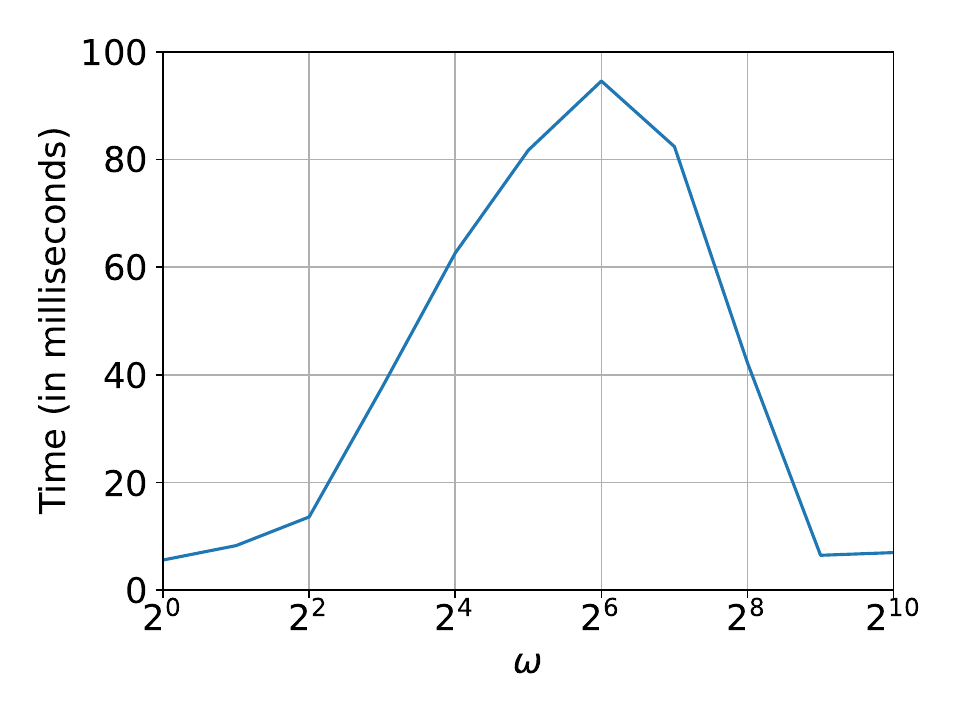}
\hfil
\includegraphics[width=.32\textwidth]{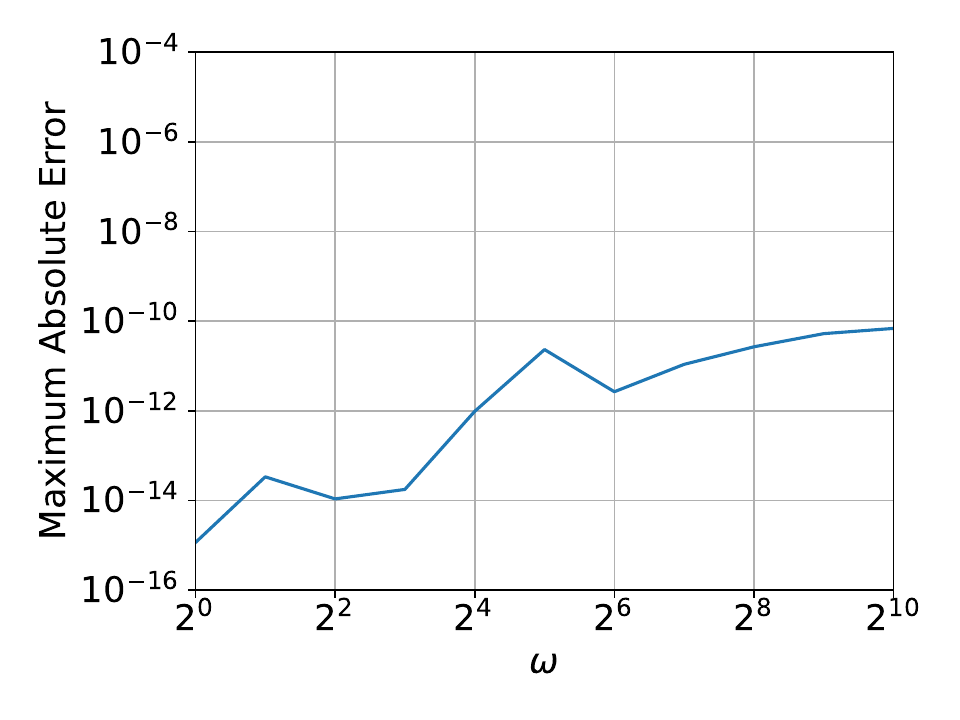}
\hfil
\includegraphics[width=.32\textwidth]{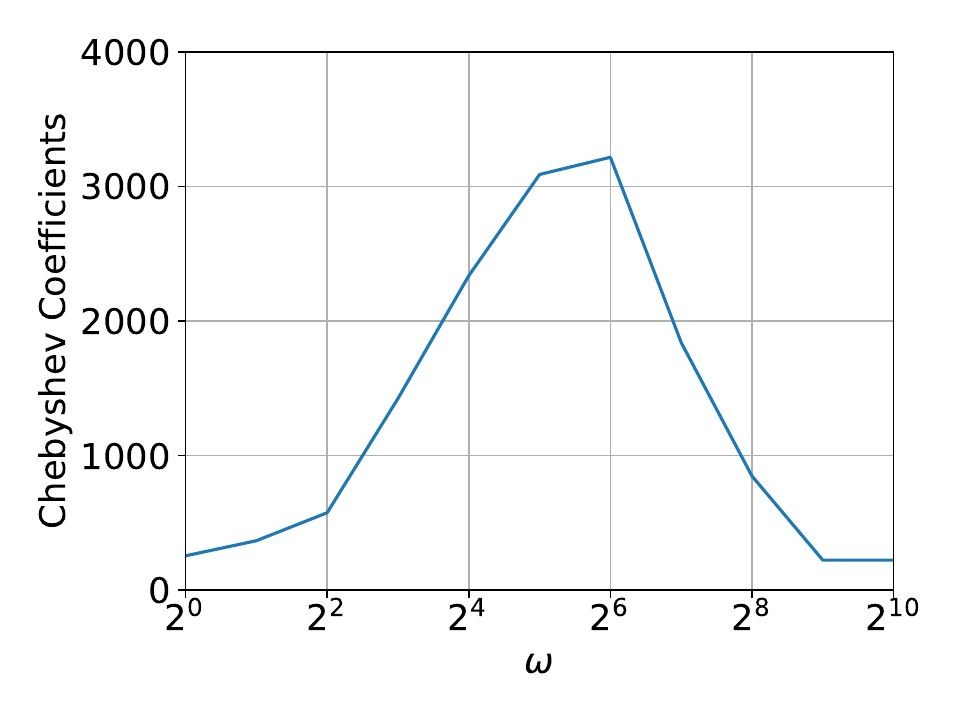}
\hfil

\hfil
\includegraphics[width=.32\textwidth]{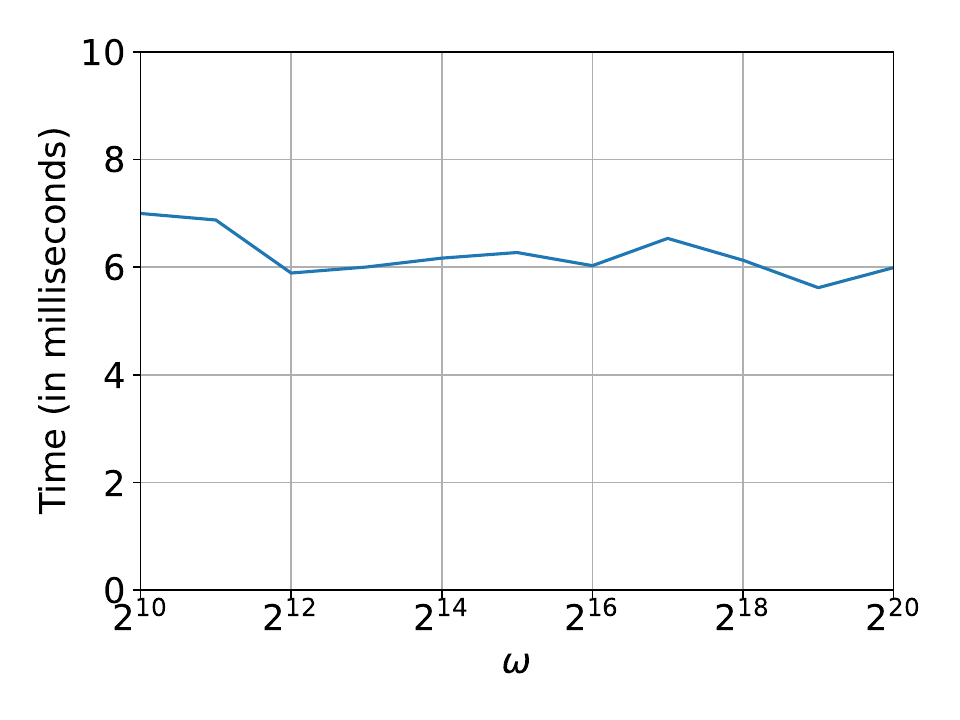}
\hfil
\includegraphics[width=.32\textwidth]{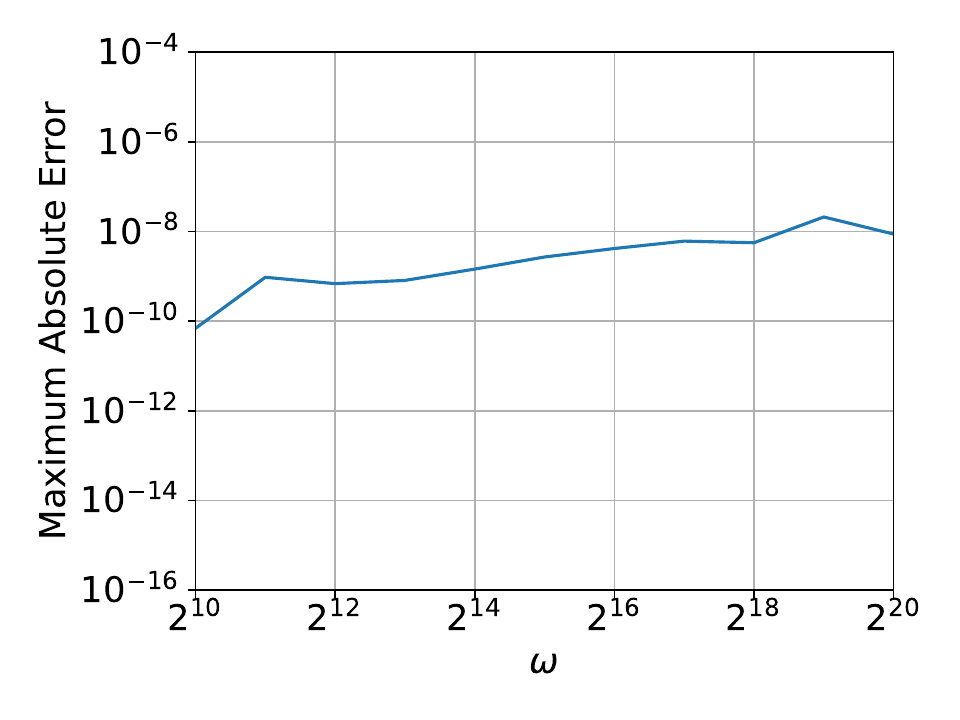}
\hfil
\includegraphics[width=.32\textwidth]{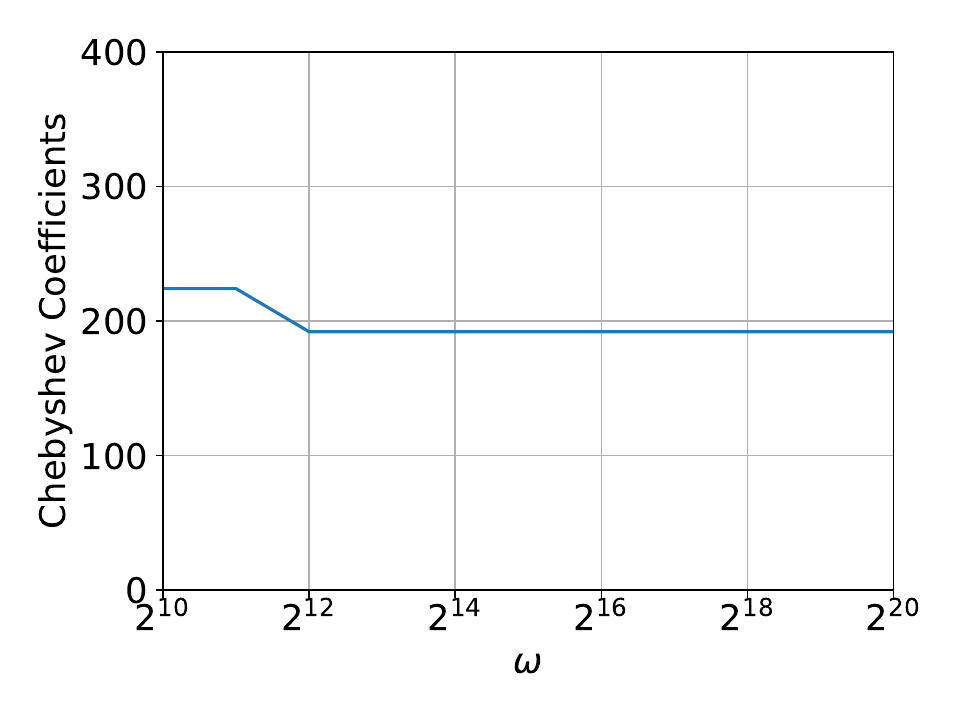}
\hfil

\caption{The results of the experiments of Subsection~\ref{section:experiments:4}.
The plot at top left gives the running time of the algorithm of this paper
in the low-frequency regime.
The top-middle plot gives reports the absolute error in the solution of the initial value problem for (\ref{experiments4:1})
in the low-frequency regime.
The plot at top right shows  the total number  of piecewise Chebyshev coefficients required to represent the slowly-varying phase functions,
 again in the low-frequency regime.
The plots on the bottom row provide the same information, but in the high-frequency regime.
}

\label{experiments4:figure1}
\end{figure}

\end{subsection}

\end{section}

\begin{section}{Conclusions}
\label{section:conclusion}

We have described a numerical algorithm for the solution of linear scalar ordinary differential
equations with slowly-varying coefficients whose running time is bounded independent of frequency.
It is  competitive with cutting edge methods for  second order equations, and significantly faster than 
state-of-the-art methods for higher order equations.   
The key observation underlying our algorithm is that the solutions
of scalar linear ordinary differential equations can be efficiently represented via phase functions.
One of the main differences between our algorithm and many alternative approaches is that, rather than trying to 
approximate phase functions with a series expansion or an iterative process, we construct them
by simply solving the Riccati equation numerically.

In the case of second order equations, the principles which underlie our solver have been rigorously
justified.  However, we have not yet proved the analogous results for higher
order scalar equations.  This is the subject  of ongoing work by authors.

There are a number of obvious mechanisms for accelerating our algorithm.  
Perhaps the simplest would be to replace the robust but fairly slow solver of Appendix~\ref{section:appendix}
with a faster method.  We could also exploit the symmetries possessed by the solutions of the Riccati equation.
For example, when the coefficient $q$ in the second order 
equation (\ref{introduction:two}) is real-valued, there is a pair of slowly-varying phase functions $\psi_1$ and $\psi_2$
related by complex conjugation (i.e., $\psi_1 = \overline{\psi_2}$) and it is only
necessary to construct one of these phase functions.

The authors have also developed a ``global'' variant of the algorithm
of this paper.  Rather than applying the  Levin procedure only to calculate the values of 
$r_1,\ldots,r_n$ at a single point in the solution domain, it uses it 
as the basis of an adaptive method for  calculating $r_1,\ldots,r_n$ over the entire solution domain.
This approach is generally faster than that of this paper in the event
that all of the eigenvalues $\lambda_1(t),\ldots,\lambda_n(t)$ of the coefficient
matrix for (\ref{introduction:scalarcoef}) are of large magnitude.  However, when one or more
of the eigenvalues is of small magnitude, the slowly-varying phase functions are nonunique
and the method runs into difficulties.  A preliminary discussion
of the global variant of our algorithm can be found in \cite{aubry2023}; a thorough description
of it will be given by the authors at a later data.  The authors also plan
to describe the generalization of the algorithm of \cite{BremerPhase}
to equations of the form (\ref{introduction:scalarode}) and compare it
to the method of this paper and its global variant.

It is straightforward to generalize our method
to the case of scalar differential equations which
are nondegenerate on an interval $[a,b]$ except at a finite collection of turning points.
This can be done by applying the algorithm of this paper to a collection of subintervals
of $[a,b]$.   

Finally, we note that because essentially any system of linear ordinary differential equations can be transformed
into a  scalar equation (see, for instance, \cite{PutSinger}), the algorithm of this paper can 
be used to solve a large class of systems of linear ordinary differential equations
in time bounded independent of frequency. The preprint \cite{hubremer} introduces an algorithm
based on this approach; that is, transforming a system of linear ordinary
differential equations into a scalar equation which is then solved via the algorithm of this paper.

\end{section}

\begin{section}{Acknowledgments}
JB  was supported in part by NSERC Discovery grant  RGPIN-2021-02613.
We thank Fruzsina Agocs for directing us to the version of 
the  algorithm of \cite{agocs} designed to solve Legendre's equation
used in the experiments of this paper.
\end{section}

\begin{section}{Data availability statement}
The datasets generated during and/or analysed during the current study are available from the corresponding author on reasonable request.
\end{section}

\bibliographystyle{acm}
\bibliography{scalar.bib}

\appendix
\begin{section}{An adaptive spectral solver for ordinary differential equations}
\label{section:appendix}

In this appendix, we detail a standard  adaptive  spectral method for solving ordinary
differential equations.  It is used  by the algorithm of this paper, and also to calculate reference
solutions in our numerical experiments. We describe its operation in the 
case of the initial value problem
\begin{equation}
\left\{
\begin{aligned}
\bm{y}'(t) &= F(t,\bm{y}(t)), \ \ \ a < t < b,\\
\bm{y}(a) &= \bm{v}
\end{aligned}
\right.
\label{algorithm:system}
\end{equation}
where $F:\mathbb{R}^{n+1} \to \mathbb{C}^n$ is smooth and $\bm{v} \in \mathbb{C}^n$.
However, the solver can be easily modified to produce a solution with a specified value
at any point $\eta$ in $[a,b]$.  Moreover, by running the solver multiple times,  a basis in the
space of solutions of a system of differential equations can be constructed
and used to  solve boundary value problems as well.

The solver takes as input a positive integer $k$, a tolerance parameter $\epsilon$, an interval $(a,b)$, 
the vector $\bm{v}$ and a 
 subroutine for evaluating the function $F$.  It outputs $n$ piecewise $(k-1)^{st}$ order Chebyshev expansions,
one for each of the components $y_i(t)$ of the solution $\bm{y}$ of (\ref{algorithm:system}).

The solver maintains two lists of subintervals of $(a,b)$: one consisting of what we term ``accepted subintervals''
and the other of subintervals which have yet to be processed.  A subinterval is accepted if the solution
is deemed to be adequately represented by a $(k-1)^{st}$ order Chebyshev expansion on that subinterval.
Initially, the list of accepted subintervals is empty and the list of 
subintervals to process contains the single interval $(a,b)$.
It then operates as follows until the list of subintervals to process is empty:
\begin{enumerate}

\item
Find, in the list of subinterval to process, the interval $(c,d)$ such that
$c$ is as small as possible and remove this subinterval from the list.

\item
Solve the initial value problem
\begin{equation}
\left\{
\begin{aligned}
\bm{u}'(t) &= F(t,\bm{u}(t)), \ \ \ c< t < d,\\
\bm{u}(c) &= \bm{w}
\end{aligned}
\right.
\label{algorithm:ivp2}
\end{equation}
If $(c,d) = (a,b)$, then we take $\bm{w}=\bm{v}$.  Otherwise,
the value of the solution at the point $c$ has already been approximated, and we use that estimate
for $\bm{w}$ in (\ref{algorithm:ivp2}).

If the problem is linear, a straightforward Chebyshev integral equation method is used to solve (\ref{algorithm:ivp2}).  
Otherwise,  the trapezoidal method is first used to produce an initial
approximation $\bm{y_0}$ of the solution and then Newton's method is applied to refine it.
The linearized problems are solved using a Chebyshev integral equation method.

In any event, the result is a set of $(k-1)^{st}$ order Chebyshev expansions 
\begin{equation}
u_i(t)  \approx \sum_{j=0}^{k-1} \lambda_{ij}\ T_j\left(\frac{2}{d-c} t + \frac{c+d}{c-d}\right),\ \ \ i=1,\ldots,n,
\label{algorithm:exps}
\end{equation}
which purportedly approximate  the components $u_1,\ldots,u_n$ of the solution of (\ref{algorithm:ivp2}).

\item
Compute the quantities
\begin{equation}
\frac{\sqrt{\sum_{j=k-2}^{k-1} \left|\lambda_{ij}\right|^2}}{\sqrt{\sum_{j=0}^{k-1} \left|\lambda_{ij}\right|^2}}, \ \ \ i=1,\ldots,n,
\end{equation}
where the $\lambda_{ij}$ are the coefficients in the expansions (\ref{algorithm:exps}).
If any of the resulting values is  larger than $\epsilon$,
then we split the subinterval into two halves $\left(c,\frac{c+d}{2}\right)$ and 
$\left(\frac{c+d}{2},d\right)$ and place them on the list of subintervals to process.  Otherwise, we place the subinterval
$(c,d)$ on the list of accepted subintervals.

\end{enumerate} 

At the conclusion of this procedure,  we have $(k-1)^{st}$ order piecewise Chebyshev expansions
for each component of the solution, with the list of accepted subintervals determining the
partition of $[a,b]$.

\end{section}

\end{document}